%
%
 

\documentstyle{amsppt} 
\TagsOnLeft
\magnification = \magstep1 
\overfullrule0pt
\hcorrection{.25in}
\advance\vsize-.75in
\def \n {\noindent}
\def \s {\smallskip}
\def \m {\medskip}
\def \b {\bigskip}
\def \a {\alpha}
 
\def \be {\beta}

\def \del {\delta}

\def\eqdef{\buildrel\text {def} \over =}  
\def\fsl{\frak s\frak l_2} 
\def\fst{\frak s\frak l_3}
\def\fss{\frak s\frak l(1,1)}

\def \ga {\gamma}  
\def \g {{\frak g}} 
\def \h {{\frak h}}
\def \k {\kappa}

\def \Hom {{\text {Hom}}}
\def \l {\lambda}
\def \la {\langle}
\def \o {\overline}

\def \ra {\rangle}

\def \u {\underline}
\def \w {\omega} 
\def \z {\zeta}
\def \C {{\text {\bf C}}}
\def \K {{\text {\bf K}}} 
\def \Z {{\text {\bf Z}}}

\topmatter
\title Down-up Algebras
\endtitle
\author Georgia Benkart \\
Tom Roby
\endauthor
\thanks The first author gratefully acknowledges 
support from  National Science 
Foundation Grant \#{}DMS--9622447 and from
the Mathematical Sciences Research Institute, Berkeley. Research at MSRI
was supported in part by NSF grant 9701755.  
\newline
The second author gratefully acknowledges 
 support from National Science Foundation Grant \#{}DMS--9353149 and the
University of Wisconsin, Madison.  
\endthanks

\leftheadtext{ GEORGIA BENKART, TOM ROBY}
\rightheadtext {DOWN-UP ALGEBRAS}
\abstract
The algebra generated by the down and up operators
 on a differential partially ordered set (poset) encodes
 essential enumerative and structural properties of the
 poset.  Motivated by the algebras generated
by the down and up operators on posets, we introduce here a family of
infinite-dimensional associative algebras called
down-up algebras.   We show that down-up algebras
exhibit many of the important features of the universal enveloping
algebra  $U(\fsl)$ of the Lie algebra $\fsl$ including a
Poincar\'e-Birkhoff-Witt type basis and a well-behaved representation theory. 
We investigate the structure and representations of 
down-up algebras and focus especially on 
Verma modules, highest weight representations, and category
$\Cal O$ modules for them. 
We calculate the exact expressions for all the weights, since that information has
proven to be particularly useful in determining structural results about posets. 
\endabstract
\endtopmatter 
\b
\head  {\S 1. Introduction}  \endhead
\b
\m
\subhead {Combinatorial source of down-up algebras} \endsubhead
\b
Assume $P$ is a partially ordered
set (poset), and let 
$\C P$  denote the complex vector space whose basis is the set $P$.  
For many posets there are two  well-defined  linear transformations on $\C P$  
 coming from the order relation on $P$,
$$ d(y)=\sum_{x \prec y} x \quad \quad  \text {and} 
\quad \quad u(y)= \sum_{y \prec z} z,$$   
 
\n where $\prec$ denotes the covering relation in the poset. 
Thus, $d(y)$ is the sum of the elements which $y$ covers, and
$u(y)$ is the sum of the elements which cover $y$.   
(These are the ``down'' and ``up'' operators.) 
Young's lattice of all partitions of all
nonnegative integers provides an important example of such
a poset.  If
$\mu = (\mu_1, \mu_2, \dots)$ and $\nu = (\nu_1, \nu_2, \dots)$
are partitions (with parts arranged so $\mu_1  \geq \mu_2 \geq \cdots$
and $\nu_1 \geq \nu_2 \geq \cdots)$, then $\mu \leq \nu$ if $\mu_i \leq \nu_i$
for each $i$.   The partition $\nu$ covers $\mu$ if $\mu \leq \nu$
and $\sum_i \nu_i = 1 +\sum_i \mu_i$.   
In ([St1], [St2]), Stanley found that many
of the interesting enumerative and structural properties of Young's
lattice could be deduced from the following relation that the
up and down 
operators satisfy,

$$ du-ud = I,$$

\n where $I$ is the identity transformation
on $\C P$. 
Young's lattice of partitions is not the only
example, and Stanley considered more general posets 
 which satisfy the relation $du-ud = rI$ for some fixed positive
integer $r$.  Such posets afford a representation of the
Weyl algebra generated by $u$ and $d/r$.  Since
the Weyl algebra also can be realized
as differential operators $y \mapsto d/dx$ and $x \mapsto x$ (multiplication
by $x$)  on $\C[x]$, Stanley referred to the posets satisfying $du-ud = rI$
as ``$r$-differential'' or simply ``differential'' when $r = 1$.  
Fomin [F] independently defined
essentially the same class of posets for $r = 1$  calling them ``Y-graphs'', the
terminology inspired by  Young's lattice.   

\b In his study of
{\it uniform posets} [T], Terwilliger considered finite ranked posets $P$
whose down and up operators satisfy the
following relation 

$$ d_{i}d_{i+1}u_{i}=\alpha _{i}d_{i}u_{i-1}d_{i}+\beta
_{i}u_{i-2}d_{i-1}d_{i} +\gamma _{i}d_{i}, $$

\n where $d_{i}$ and $u_{i}$ denote the restriction of $d$ and $u$ to the
elements of rank $i$. (There is an analogous second relation,

$$ d_{i+1}u_{i}u_{i-1}=\alpha _{i}u_{i-1}d_{i}u_{i-1}+\beta
_{i}u_{i-1}u_{i-2}d_{i-1} +\gamma _{i}u_{i-1}, $$
 
\n which holds automatically in this case because $d_{i+1}$ and $u_i$ are
adjoint operators relative to a certain bilinear form.)
The matrix algebra over the reals generated by these
operators together with the idempotent projection maps onto each rank,
(referred to as the ``incidence algebra'' in [T]), is finite-dimensional
and in [T] its finite-dimensional simple representations are
determined.  Terwilliger also discusses many families of examples including
the four presented below.  

\b   
In many classical cases the constants in the above relations do not
depend on the rank of the poset.   A particular instance of this
provides a $q$-analogue
of the notion of differential poset.  A 
partially ordered set whose down and 
up operators satisfy 

$$
\aligned
d^{2}u & = q(q+1)dud -q^{3}ud^{2}+rd \\
du^2 & = q(q+1)udu - q^3 u^2 d + r u \\
\endaligned \tag 1.1$$

\n where $q$ and $r$ are fixed complex numbers is 
said to be  ``$(q,r)$-differential''.  
 Examples of $(q,r)$-differential posets 
include the posets of alternating forms,  
quadratic forms, or Hermitian forms over a finite field,  as well as Hemmeter's
poset: 
\b  
\roster
\item  {[$Alt_{q}(n)$, $(n \geq 2)$]}\ {\it The poset of alternating 
bilinear forms}
\s \n  Assume $H$ is an $n$-dimensional vector space over the field
$GF(q)$ of $q$ elements.   Consider the set of pairs $P=\{(U,f)
\mid U$ is a subspace of $H$ and $f$ is an alternating bilinear form on $U
\}$ with the following ordering: $(U,f)\leq (V,g)$ if $U$ is a subspace
of $V$ and $g|_{U}=f$.  Then $P$ is a $(q,r)$-differential poset with
$r=-q^n(q+1)$.  
\m
\item  {[$Quad_{q}(n)$, $(n \geq 2)$]}\ {\it The poset of quadratic forms}
\s   
\n The poset is defined the same way as in the first example but with 
``quadratic form'' replacing ``alternating bilinear form''.  In this case
$P$ is
\hskip 0pt plus 20pt\penalty-999
$(q,\nobreak -q^{n+1}(q+\nobreak 1))$-differential.  
\m
\item  {[$Her_{q}(n)$, $(n \geq 2)$]}\ {\it The poset of Hermitian forms} 
\s    
\n In this example, $H$ is taken to be an $n$-dimensional vector space over 
$GF(q^{2}),$  and the bilinear forms are ``Hermitian''.  The result is a poset
which is
$(q^{2},-q^{2n+1}(q^{2}+1))$-differential.    
\m
\item  {[$Hem_{q}(n)$, ($q$ odd, $n\geq 2$)]}\ {\it Hemmeter's poset}
\s    
\n Suppose $X^{+}$ and $X^{-}$ are two copies of the graph of the dual polar
space $X=C_{n-1}(q)$.  Let $Y$ be the bipartite graph with vertex set
$X^{+}\cup X^{-}$ (forming the bipartition).   Two vertices $x^{+}\in
X^{+}$ and $y^{-}\in X^{-}$ are adjacent in $Y$ if and only if  $x=y$ or $x$ 
and $y$ are adjacent in $X$.   The graph $Y$ is a
distance regular graph recently discovered by Hemmeter.   Now set $P=Y$ and fix
$u_{0}\in Y$.  Declare
$u\leq v$ if $\del (u_{0},u)+\del (u,v)=\del (u_{0},v)$,
where $\del(u,v)$ is the distance between $u$ and $v$.  (The choice of $u_{0}$
is irrelevant since $Y$ is vertex transitive.)   
Then like the poset $Alt_{q}(n)$,
$P$ is $(q,-q^n(q+1))$-differential.  
\endroster
\m
\b 
\subhead {Definition of down-up algebras}\endsubhead
\b
In this paper we study certain infinite-dimensional associative
algebras whose generators satisfy relations more general than (1.1).  
Although
the original motivation for our investigations came from posets,
we make no assumption about the existence of posets whose
down and up operators satisfy our relations.

\b 
We say a unital associative
algebra 
$A = A(\a,\be,\ga)$ over the complex numbers $\C$
with generators $u,d$ and defining relations
\b 
\item {}{(R1)} $d^2u = \a dud + \be ud^2 + \ga d,$
\m 
\item {} {(R2)} $du^2 = \a udu + \be u^2 d + \ga u,$ 
\b 
\n 
where $\a,\be,\ga$ are fixed but arbitrary elements of $\C$,
is a {\bf down-up algebra}.  
\b
It is easy to see that when $\ga \neq 0$ the down-up algebra $A(\a,\be,\ga)$
is isomorphic to $A(\a,\be,1)$ by the map  $d \mapsto d'$, $u \mapsto \ga u'$.
Therefore, it would suffice to treat just two cases $\ga = 0,1$.  But rather
than divide all our considerations into these cases, it is convenient to retain the notation
$\ga$ throughout.   
\b
\subhead {Examples of down-up algebras}\endsubhead
\b
If
$d$ and
$u$ are the down and up operators of
a $(q,r)$-differential poset, and $B$ is
the associative algebra they generate,
then $B$ satisfies (R1) and (R2) with 
$\a = q(q+1)$, $\be = -q^{3}$, and $\ga =  r$.  Thus, $B$ is a homomorphic image
of the algebra $A(\a,\be,\ga)$ with these parameters.
In an $r$-differential poset, the relation $du - ud = rI$
is satisfied.  Multiplying on the left by $d$ and on the right by $d$ and
adding the resulting equations, we get the relation
$d^2 u - ud^2 = 2rd$ of a $(-1,2r)$-differential poset.  Thus,
the Weyl algebra is a homomorphic image (by the ideal generated
by $du-ud - r1$)  of the algebra $A(0,1,2r)$. More generally,
the $q$-Weyl algebra is a homomorphic image of the algebra 
$A(0,q^2,q+1)$ (and also of the algebra $A(q-1,q,2)$) 
by the ideal generated by $du - qud - 1$.  
The skew polynomial
algebra $\C_q[d,u]$, or quantum plane (see [M]), is the associative algebra with
generators $d,u$ which satisfy the relation  $du = qud$.  Therefore,
$\C_q[d,u]$ is a homomorphic image (by the ideal generated by $du-qud$)
of the algebra $A(2q,-q^2,0)$ (also of the algebra $A(aq,(1-a)q^2,0)$ for any $a \in \C$).  

\b 
Consider a 3-dimensional Lie algebra $L$ over $\C$ having a basis $x,y,[x,y]$ 
such that $[x[x,y]] = \ga x$ and  $[[x,y],y] = \ga y$.  
In the universal enveloping
algebra $U(L)$ of $L$, the relations above become

$$\aligned
& x^2 y -2 xyx + yx^2 = \ga x \\
& xy^2 -2yxy +  y^2x = \ga y. \\
\endaligned$$

\n Thus, $U(L)$ is a homomorphic image of the down-up algebra $A(2,-1,\ga)$
via the mapping $\phi: A(2,-1,\ga) \rightarrow U(L)$ with $\phi: d \mapsto x$, $\phi: u
\mapsto y$. In fact, $U(L)$ is isomorphic to $A(2,-1,\ga)$.  The
mapping $\psi: L \rightarrow A(2,-1,\ga)$ with $\psi: x \mapsto d$, $
\psi: y \mapsto u$,
and $\psi: [x,y] \mapsto du - ud$ extends, by the universal property of $U(L)$,
to an algebra homomorphism $\psi: U(L) \rightarrow A(2,-1,\ga)$.  The 
two homomorphisms $\phi$ and $\psi$ are inverses of each other.  
Every 3-dimensional Lie algebra $L$ over $\C$ with $L = [L,L]$ has such a basis (see
[J, pp. 13-14]), 
so its enveloping algebra is a down-up algebra.  
\b
There are several noteworthy special cases. The Lie algebra
$\fsl$ of $2 \times 2$ complex matrices of trace zero has a standard
basis $e,f,h$, which satisfies $[e,f] = h, \; [h,e] = 2e$, and $[h,f] = -2 f$. From
this we see that $U(\fsl) \cong A(2,-1,-2)$.  The Heisenberg Lie algebra $\frak H$ has
a basis $x,y,z$ where $[x,y] = z$, and $[z,{\frak H}] = 0$. Thus, $U({\frak H}) \cong
A(2,-1,0)$.  
\b
In [Sm], S.P. Smith investigated a class of associative algebras that have
many of the properties of the enveloping algebra $U(\fsl)$ but also
many new, interesting aspects.  These algebras have a presentation by
generators
$a,b,h$ and relations $[h,a] = a$, $[h,b] = -b$ and $ab - ba = f(h)$,
where $f(h)$ is a polynomial in $h$. In the special situation that
$\deg(f) \leq 1$,  such an algebra is a homomorphic image of a down-up algebra
$A(2,-1,\gamma)$ for some $\gamma$.   Algebras related to Smith's appear
in Hodges' work [H]  on primitive factors of $U(\fsl)$ and type-A Kleinian
singularities. 
\b
Woronowicz [Wo] introduced a family of algebras whose
presentation is equivalent to having three generators
$x,y,z$ which are assumed to satisfy the relations 

$$x z - \z^{4} z x = (1+\z^{2})x \quad z y - \z^{4}yz = (1+\z^{2})y
\quad xy - \z^{2}y x = \z z$$

\n for some scalar $\z \neq \pm 1, 0$.  These algebras are down-up algebras with $\a = 
\z^{2}(1+\z^2)$, \ $\be = -\z^{6}$, and  $\ga = \z(1+\z^{2})$.
Woronowicz's algebras fit into a more general construction studied by
Jordan in [Jo].    
\b
The Lie superalgebra $S = \fss = S_{\o 0} \oplus S_{\o 1}$ of $2 \times 2$ complex
matrices, 
$y = \left (\matrix  y_1 & y_2 \\  y_3 & y_4  \endmatrix \right )$
with supertrace $y_1-y_4 = 0$ and with multiplication
given by the supercommutator $[x,y] = xy-(-1)^{ab}yx$ for $x\in S_{\o a},\;y\in S_{\o
b}$  has a presentation  by generators $e,f$ (which are odd and can be identified
with the matrix units $e = e_{1,2}$, $f = e_{2,1}$) and relations
$[e,[e,f]] = 0$, $[[e,f],f] = 0$, $[e,e] = 0$, $[f,f] = 0$.
Its universal enveloping algebra $U(\fss)$ has generators $e,f$ and
relations $e^2 f - f e^2 = 0$, $e f^2 - f^2 e = 0$, $e^2 = 0$, $f^2 = 0$.
Thus, $U(\fss)$ is a homomorphic image of the down-up algebra
$A(0,1,0)$ by the ideal generated by the elements $e^2$ and $f^2$, which are
central in $A(0,1,0)$. 
\b
\m
\subhead {Connections with Witten's deformation of $U(\fsl)$}
\endsubhead
\b
To provide an explanation of the existence of quantum groups,
Witten ([W1], [W2]) introduced a 7-parameter deformation of the
universal enveloping algebra $U(\fsl)$.
Witten's deformation is a unital associative algebra
over a field $\K$ (which is algebraically closed of
characteristic zero and which could be
assumed to be $\C$) and depends on a 7-tuple $\u \xi =
(\xi_1, \dots, \xi_7)$ of elements of $\K$.  It has
a presentation by generators $x,y,z$ and defining
relations

$$\gather 
 xz - \xi_1 z x = \xi_2 x  \qquad \qquad  zy - \xi_3 yz = \xi_4 y  \\
yx - \xi_5 xy = \xi_6 z^2 + \xi_7 z. \\
\endgather$$

\n We denote the resulting algebra by $ {\frak W}(\u \xi)  $.
It can be argued that the following holds:
\b
\proclaim{Theorem 1.2} ([B, Thm. 2.6]) A Witten deformation algebra ${\frak W}(\u
\xi) 
$ with

$$\xi_6 = 0, \;\; \xi_5\xi_7 \neq 0, \; \xi_1 = \xi_3, \; \;\text{and}
\;\;\xi_2 = \xi_4 \tag 1.3$$

\n is isomorphic to the down-up algebra $A(\a,\be,\ga)$
with $\a,\be,\ga$ given by   
 
$$\a = \frac{1+\xi_1 \xi_5}{\xi_5}, \quad \quad 
\be = - \frac{\xi_1}{\xi_5}, \quad \quad \gamma = -\frac{\xi_2
\xi_7}{\xi_5}. \tag 1.4$$ 

\n Conversely, any
down-up algebra $A(\a,\be,\ga)$ with not both $\a$ and $\be$
equal to 0 is isomorphic to a Witten deformation algebra ${\frak W}(\u \xi)  $
whose parameters satisfy (1.3).   
\endproclaim
\b
A deformation algebra ${\frak W}(\u \xi)$
has a filtration,  and 
Le Bruyn ([L1], [L2]) investigated the algebras ${\frak W}(\u \xi)  $ 
whose associated graded algebras are Auslander regular. 
They determine a 3-parameter family of deformation algebras which
are called {\it conformal $\fsl$ algebras} and whose defining relations
are 

$$ xz - a zx = x \quad \quad zy - a yz = y \quad \quad
yx - c xy = b z^2 + z \tag 1.5$$

\n When $c \neq 0$ and $b = 0$, the conformal $\fsl$ algebra
with defining relations given by (1.5) is isomorphic to
the down-up algebra $A(\a,\be,\ga)$ with
$\a = c^{-1}(1+ac), \be = -ac^{-1}$ and $\ga = -c^{-1}$.
If $b= c = 0$ and $a \neq 0$, then the conformal $\fsl$ algebra
is isomorphic to the down-up algebra $A(\a,\be,\ga)$
with $\a = a^{-1}$, $\be = 0$ and $\ga = -a^{-1}$.  
\b
In a recent paper Kulkarni [K1] has shown that under certain
assumptions on the parameters a Witten deformation algebra is isomorphic to
a conformal $\fsl$ algebra or to a double skew polynomial extension.
The precise statement of the result is
\b
\proclaim {Theorem 1.6} ([K1, Thm. 3.0.3]) If  $\xi_1\xi_3\xi_5 \xi_2 \neq 0$ or
$\xi_1
\xi_3
\xi_5 \xi_4 \neq 0$, then ${\frak W}(\u \xi)$ is isomorphic to one of the following
algebras:
\s
\item {(a)} A conformal $\fsl$ algebra with generators $x,y,z$ and relations
given by (1.5) for some $a,b,c \in \K$.  
\s
\item {(b)} A double skew polynomial extension (that is,
a skew polynomial extension of a skew polynomial ring) whose generators satisfy
\s \item{} \quad {(i)} $xz - zx = x, \quad zy - yz = \z y, \quad yx - \eta xy = 0$ \quad or
\s \item{} \quad {(ii)} $xw = \theta wx$, \quad $wy = \kappa yw$, \quad $yx = \lambda xy$
\s \item{} for parameters $\z,\eta,\theta,\kappa,\lambda \in \K.$ \endproclaim 

\b  
Kulkarni [K1] studies the simple representations of
the conformal $sl_2$ algebras and of the skew polynomial
algebras in (b) by making essential use of the
observation that the conformal $\fsl$ algebras of (1.5) can be realized as {\it
hyperbolic rings}.  Kulkarni applies results of Rosenberg [R] on left ideals of
the left spectrum of hyperbolic rings to
determine the maximal left ideals for the conformal $sl_2$ algebras. Further
applications of results of [R] give the left spectrum of the double skew
polynomial extensions in (b) of the theorem.    
\m   
\b 
In this work we investigate the structure and representations of the
down-up algebras $A(\a,\be,\ga)$.  Remarkably down-up algebras have
many of the important features of universal enveloping algebras including a 
Poincar\'e-Birkhoff-Witt basis theorem and a well-behaved representation theory.
Adopting a Lie theoretic approach, we explicitly construct
``highest weight'' and ``lowest weight'' representations of the algebras
$A(\a,\be,\ga)$ and determine when they are simple.  We calculate 
the exact expressions for all the weights since that information has proven to be
especially useful in determining structural results about posets.  
(In particular, for differential posets Stanley [St1] showed that the spectrum of
$ud$ is related to the eigenvalues of the adjacency matrix of
certain finite graphs associated with the posets. This gave a new class
of graphs whose spectra could be explicitly computed.  Many interesting
properties associated to random walks on the graph are
related to its spectrum.) 
\b
In the spirit
of the work of [BGG] we investigate two categories of modules ($\Cal O$ and $\Cal O'$)
for down-up algebras whose
simple objects include the highest weight simple modules.   Corresponding
to each complex number $\l$, there is a universal 
infinite-dimensional highest weight module $V(\l)$, 
the so-called Verma module.  Each of the weight spaces of $V(\l)$ is 1-dimensional
or infinite-dimensional depending on the choice of the parameters
$\a,\be,\ga,\l$.  (Roots of unity play a critical role
in deciding which occurs.)   Its weights are determined by
the recurrence relation

$$\l_n = \a \l_{n-1} + \be \l_{n-2} + \ga, \quad \quad \l_0 = \l,\quad \l_{-1} = 0.$$

\n In particular, for the algebra $A(1,1,0)$, the weights of
$V(\l)$  are given 
by the Fibonacci sequence $\l,\l,2\l, 3\l,5\l,8\l, \dots$ for
each $\l$. 
Besides the examples coming from $(q,r)$-differential posets 
and universal enveloping
algebras,  certain subalgebras of quantized enveloping
algebras (quantum groups) are down-up algebras.  The highest
weight representations that we construct for these subalgebras
are not the standard ones. 
\b
Our investigations arose from conversations one of the authors (TR)
had with Paul Terwilliger.  We would like to express our appreciation
to him for many helpful suggestions over the course of this
project. 
\b
\b
\head {\S 2. Representations of down-up algebras} \endhead
\b 
As before, let 
$A = A(\a,\be,\ga)$ be a down-up algebra over $\C$
with generators $u,d$ and defining relations
\m
\item {}{(R1)} $d^2u = \a dud + \be ud^2 + \ga d,$ 
\s
\item {}{(R2)} $du^2 = \a udu + \be u^2 d +
\ga u$, 
\m
\n where $\a,\be,\ga$ are fixed but arbitrary elements of
$\C$. 
\b
\m
\subhead Highest weight modules \endsubhead
\b
Our first objective is to construct
``highest weight''  representations of down-up algebras.  
Involved in the construction is 
a certain sequence of elements of $\C$. Starting
with $\l_{-1} = 0$ and an arbitrary $\l_0 = \l \in \C$, 
define $\l_n$
for $n \geq 1$ inductively by the following
recurrence relation: 

$$\lambda_n = \a \lambda_{n-1} + \be \l_{n-2} + \ga. \tag 2.1$$

\m
\proclaim {Proposition 2.2} Set $\l_{-1} = 0$. Let
$\l_0 = \l \in \C$ be arbitrary, and
define $\l_n$ for $n \geq 1$ as in (2.1).   Let $V(\l)$ be the
$\C$-vector space having basis $\{v_n \mid n = 0,1,2, \dots\}$.  Define

$$\aligned
d \cdot v_n & = \l_{n-1} v_{n-1}, \quad n \geq 1, \quad \text {and}
\quad d
\cdot v_0 = 0 \\
u \cdot v_n & = v_{n+1}. \\
\endaligned \tag 2.3$$

\n Then this action extends to an $A(\a,\be,\gamma)$-module
action on $V(\l)$.  \endproclaim
\b
{\bf Proof.} The free associative algebra 
$\C\langle u,d\rangle$ with $1$ over
$\C$ generated by $u,d$ acts on $V(\l)$ by extending
(2.3) making $V(\l)$ into an $\C\langle u,d\rangle$-module. 
 To show there is an induced
action of $A$, it suffices to verify that

$$\aligned 
&(d^2 u - \a dud - \be ud^2 - \ga d)v_n = 0 \\
&(du^2 -\a udu - \be u^2 d - \ga u)v_n = 0 \\
\endaligned$$

\n for all $n \geq 0$.   But this is easy to see since 

$$(d^2 u - \a dud - \be ud^2 - \ga d)v_n =
\l_{n-1}(\l_n - \a \l_{n-1} - \be \l_{n-2} - \ga)v_{n-1}, $$

\n and

$$(du^2 -\a udu - \be u^2 d - \ga u)v_n = (\l_{n+1}
-\a \l_n - \be \l_{n-1}-\ga)v_{n+1} = 0. \qed$$
\b
\proclaim {Proposition 2.4}
 \item {}{(a)} The $A(\a,\be,\ga)$-module $V(\l)$ is simple
if and only if $\l_n \neq 0$ for any $n$. 
\m
\item{}{(b)} If $m$ is minimal with the property that
$\l_m = 0$, then $M(\l) =$ span$_\C\{v_j \mid j \geq m+1\}$ is a maximal
submodule of $V(\l)$.   
\m
\item{}{(c)} Suppose $N$ is a submodule of $V(\l)$ such that
$N \subseteq$ span$_\C\{v_j \mid j \geq 1\}$.  Then $N \subseteq M(\l)$. 

 \endproclaim 
\b
{\bf Proof.} Assume $N$ is a nonzero submodule of $V(\l)$, and  
suppose $x = \sum_{k = 0}^n a_k v_k \in N$, where $a_n \neq 0$.
Then 

$$d^{n} \cdot x = \l_{n-1} \l_{n-2} \cdots \l_1 \l_0 a_n v_0 \in N.$$

\n If the coefficient is nonzero, then $v_0 \in N$, and since
$v_0$ generates $V(\l)$, we have $N = V(\l)$ in this case.
 
If $\l_r = 0$, then 
span$_\C\{v_j \mid j \geq r+1\}$ is a proper 
submodule of $V(\l)$ by (2.3).  Thus, $V(\l)$ is simple
if and only if $\l_r \neq 0$ for any $r$.  

Let $m$ be the minimal
value such that $\l_m = 0$.  Then $M(\l)
=$ span$_\C\{v_j \mid j \geq m+1\}$ is a
submodule of $V(\l)$.   Suppose $N$ is
a submodule of $V(\l)$ properly containing $M(\l)$.  
We may assume that some $y = \sum_{j = 0}^m b_j v_j \neq 0$ belongs to
$N$. If
$n \leq m$ is maximal so $b_n \neq 0$, then $d^{n} \cdot y
= \l_{n-1} \l_{n-2} \cdots \l_1 \l_0 b_n v_0 \in N$. 
By the minimality of $m$ this implies $v_0 \in N$ and
hence that $N = V(\l)$.  Thus, $M(\l)$ is a maximal submodule
of $V(\l)$.  

Now assume $N$ is a submodule of $V(\l)$ contained in $V_1 =$ 
span$_\C\{v_j \mid j \geq 1\}$.  The submodule $N + M(\l)$ is
contained in $V_1$, so it is a proper submodule of $V(\l)$.
By the maximality of $M(\l)$, we must have $N+M(\l) = M(\l)$,
which implies $N \subseteq M(\l)$.    \qed
\b
\n {\bf Definition 2.5.} \quad We will set $M(\l) = (0)$ if $V(\l)$ is simple
and otherwise assume $M(\l)$ is defined by Proposition 2.4 (b) using the minimal
$m$ such that $\l_m = 0$. Note that when $M(\l) \neq (0)$, then
$M(\l) \cong V(\l_{m+1})$.    
\b
\n {\bf Example 2.6.} Observe for any down-up algebra $A = A(\a,\be,\ga)$ that
$M(0) = \text{span}_\C \{v_j \mid j \geq 1\}$, so that $V(0)/M(0)$ is
one-dimensional.    There is an algebra homomorphism $\epsilon: A \rightarrow \C$
with $\epsilon(u) = 0 = \epsilon(d)$ and $\epsilon(1) = 1$.  This gives
rise to a one-dimensional $A$-module, $\C v$, such that $a \cdot v = \epsilon(a)v$
for all $a \in A$, and $\C v \cong V(0)/M(0)$.    
\b
\proclaim {Definition 2.7}  A module $V$ for $A = A(\a,\be,\ga)$
is said to be a {\bf highest weight module of weight $\l$}
if $V$ has a vector $y_0$ such that $d \cdot y_0 = 0$,
$du \cdot y_0 = \l y_0$, and $V = Ay_0$.   The vector
$y_0$ is a {\bf maximal vector} or {\bf highest weight vector} of
$V$.  \endproclaim
\b
The highest weight module $V(\l)$ is said to be the {\bf Verma module}
of highest weight $\l$.   It is so named because it shares the same
universal property as Verma modules for finite-dimensional
semisimple complex Lie algebras: 
\b
\proclaim {Proposition 2.8}  If $V$ is a highest weight module
of weight $\l$ of $A = A(\a,\be,\ga)$, then $V$ is a homomorphic image of
$V(\l)$.   \endproclaim 
\b
{\bf Proof.} Let $y_0$ be a maximal vector of $V$ of weight
$\l$.  Then $V = A y_0$.  The mapping $V(\l) \rightarrow V$
with $a \cdot v_0 \mapsto a \cdot y_0$ for all $a \in A$ is
an $A$-module epimorphism.  \qed
\b
Suppose $V$ is any highest weight
$A$-module of weight $\l$.  Starting with a nonzero maximal vector 
$y_0 \in V$ of
weight $\l$, 
define $y_j$ inductively by $y_j = u \cdot y_{j-1}$ for $j \geq 1$.
Stop the inductive definition at the first place that a vector
$y_{m+1}$ is produced which is a linear combination of
the previous vectors.  In this case
we claim that $V = $ span$_\C\{y_j \mid j = 0,1,
\dots, m\}$.   
If no such dependency relation occurs,
then we claim that
$V =$ span$_\C\{y_j \mid j = 0,1,
\dots \},$ where the vectors $y_j$ are linearly independent.  
To see these two assertions, first we show that 

$$du \cdot y_n = \l_n y_n  \tag 2.9$$

\n where $\l_0 = \l$ and the $\l_n$'s satisfy
the recurrence relation in (2.1).   The proof is by induction on
$n$ with the inductive step being 

$$\aligned
du \cdot y_{n+1} & = du^2 \cdot y_n = \a u du \cdot y_n + \be u^2 d \cdot
y_n + \ga u \cdot y_n \\
& = \a \l_n y_{n+1} + \be u^2du \cdot y_{n-1} + \ga y_{n+1} \\
& = (\a \l_n + \be \l_{n-1} + \ga)y_{n+1} = \l_{n+1}y_{n+1}. \\
\endaligned$$

\n Observe that this implies 

$$d \cdot y_n = du \cdot y_{n-1}
= \l_{n-1} y_{n-1} \tag 2.10$$ 

\n for all $n \geq 1$.  Since $d$ and $u$ generate $A$, the relation
$y_{n+1} = u \cdot y_n$ together with (2.10) shows that 
$V =$ span$_\C\{y_j \mid j = 0,1,
\dots,m \}$ when a dependency occurs, and
$V =$ span$_\C\{y_j \mid j = 0,1,
\dots \}$ otherwise.  In both instances the $y_j$'s determine
a basis for $V$.  
\b
\m 
\subhead Weight modules \endsubhead
\b
If we multiply the relation $d^2 u - \a dud - \be u d^2 = \ga d$
on the left by $u$ and the relation $d u^2 - \a udu - \be u^2 d = \ga u$
on the right by $d$ and subtract the second from the first, the resulting
equation is 

$$0 = ud^2u - du^2d \quad \quad \text {or} \quad \quad (du)(ud) = (ud)(du).
\tag 2.11$$

\n Therefore, the  elements $du$ and $ud$ commute in $A = A(\a,\be,\ga)$.  For
any basis element $v_n \in V(\l)$, we have $du \cdot v_n = \l_n v_n$ and $ud \cdot v_n
= \l_{n-1}v_n$.  Using that with $n = 0$ and $\l \neq 0$, it is easy to see that $du$ and $ud$
are linearly independent.   Let $\h = \C du \oplus \C ud$. 
\b
We say an $A$-module $V$ is a {\bf weight module} if $V = \sum_{\nu \in \h^*} V_\nu$,
where $V_\nu = \{v \in V \mid h \cdot v = \nu(h)v$ for all $h \in \h\}$, and
the sum is over elements in the dual space $\h^*$ of $\h$.  (Necessarily
this is a direct sum.)
Any submodule of a weight module is a weight module.    
If $V_\nu \neq (0)$, then $\nu$ is a {\bf weight} and $V_\nu$ is the
corresponding {\bf weight space}. Each weight $\nu$
is determined by the pair $(\nu',\nu'')$ of complex numbers, $\nu' = \nu(du)$
and $\nu'' = \nu(ud)$, and often it is convenient to identify $\nu$ with $(\nu',\nu'')$. 
In particular, highest weight modules are weight modules in this sense.
The basis vector $v_n$ of $V(\l)$ is a weight vector whose weight
is given by the pair $(\l_n,\l_{n-1})$.  We will explore next
when two weights $(\l_k,\l_{k-1})$ and $(\l_\ell,\l_{\ell-1})$
can be equal.  For that we will need to know the values
of $\l_n$ more precisely.   
 
\b

\proclaim { Proposition 2.12} Assume  $\l_{-1} = 0$, $\l_0 = \l$,
and $\l_{n}$ for
$n \geq 1$ is given by the recurrence relation
$\l_{n}-\a\l_{n-1} -\be \l_{n-2} = \ga$.  Fix $t \in \C$ such
that
$$t^2 = \frac {\a^2 + 4\be}{4}.$$
\m
\item {}{(i)} If $\a^2 + 4 \be\neq 0$, then 

$$\l_n = c_1 r_1^n
+ c_2 r_2^n  + x_n,$$

\n where 

$$
\aligned & r_1 = \frac {\a} {2} + t, \quad
 \quad 
 r_2 = \frac {\a} {2} - t, \\
& x_n = \cases 
(1-\a-\be)^{-1}\ga \ \ \text {\quad if \quad  $\a+\be \neq 1$} \\
(2-\a)^{-1}\ga n \ \ \text {\quad if \quad $\a+\be = 1$ \quad (necessarily \quad $\a
\neq 2),$} 
\\
\endcases \\
& \text {and} \quad \quad \left (\matrix c_1 \\ c_2 \\ \endmatrix \right )
 = \frac {1}{r_2-r_1}
\left (\matrix r_2 & -1 \\ -r_1 & 1 \\ \endmatrix \right )
\left (\matrix \l-x_0 \\ \a\l+\ga-x_1 \\ \endmatrix \right ). 
\endaligned $$
 \m
\item {}{(ii)} If $\a^2 + 4 \be = 0$ and $\a \neq 0$, then 

$$\l_n = c_1 s^n + c_2 n s^n + x_n \quad \text {where}$$

$$\aligned & s = \displaystyle \frac {\a}{2} \\
& x_n = \cases 
(1-\a-\be)^{-1}\ga \ \  \text {\quad if \quad $\a+\be \neq 1$} \\
 2^{-1} n^2 \ga \ \  \text {\quad if \quad $\a + \be = 1$ 
\quad i.e. \quad if \quad  $\a = 2,\; \be = -1$,}  
\\ \endcases \\ 
& \text {and}  \quad \quad \left (\matrix c_1 \\ c_2 \\ \endmatrix \right)
 =  
\left (\matrix 1 & 0 \\ -1 & 2\a^{-1}  \\ \endmatrix \right )
\left (\matrix \l-x_0 \\ \a\l+\ga-x_1 \\ \endmatrix \right ). 
\endaligned $$
\m
\item {}{(iii)} If $\a^2 + 4 \be = 0$  and
$\a = 0$,  then $\be = 0$ and $\l_n = \ga$ for all $n \geq 1$.
\endproclaim
\b
 Note in particular, that if $\a,\be$ are
real, then it is natural to take  $t  = (1/2)\sqrt{\a^2+4\be}$.
\b 
{\bf Proof.} These assertions can be shown using standard
arguments for solving linear recurrence relations
([Br, Chap. 7]).  \qed 
\b
\proclaim{Theorem 2.13} Assume for the
algebra $A = A(\a,\be,\ga)$ that  $\l_\ell = \l_k$
and $\l_{\ell-1} = \l_{k-1}$ for some $\ell > k$.  Then
\m
\item {(a)} $A = A(\a,\be,\ga)$, where $\a^2 + 4 \be \neq 0$ and
$\a+\be \neq 1$, and  either 
\item{}{(1)} $\ga = 0$ and $\l_n = 0$
for all $n$, 

\item{}{(2)} $\ga \neq 0$, $\be = 0$, 
and $\l_n = (1-\a)^{-1}\ga$ for all $n$, 

\item{}{(3)}{ 
$$\l_n = \Big(\l- \displaystyle\frac{\ga}{1-\a-\be}\Big) r^n +
\displaystyle\frac{\ga}{1-\a-\be},$$
where $r = r_1 = (1/2)\a + t$
or $r = r_2 = (1/2)\a - t$ \ ($t^2 = (1/4)(\a^2+4\be)$),  and 
$r^{\ell-k} = 1$, or}
 
\item{}{(4)} $r_1^{\ell-k} = 1 = r_2^{\ell-k}$ and 
$$\l_n = c_1 r_1^n + c_2 r_2^n + \displaystyle\frac{\ga}{1-\a-\be},$$
where $c_1$ and $c_2$ are as in Proposition 2.12;
\m

\item{(b)} $A = A(\a,1-\a,0)$, where $\a \neq 2$, and
either 
\item{}{(1)} $\l_n = 0$ for all $n$,  or $\a = 1$ and  $\l_n = \l \neq 0$  for all $n$, 
or

\item{}{(2)}  $\l \neq 0$, $(\a-1)^{\ell-k} = 1$, and

$$\l_n = \frac {\l}{2-\a}\Big( 1 - (\a-1)^{n+1}\Big);$$
 
\m
\item {(c)} $A = A(\a,-(1/4)\a^2,\ga)$ where $\a \neq 0,2$,
and either 

\item{}{(1)} $\ga = 0$ and $\l_n = 0$ for all $n$, or

\item{}{(2)} $\ga \neq 0$, $s = (1/2)\a$ satisfies $s^{\ell-k} = 1$, and 

$$\l_n =\frac {2\ga}{(2-\a)^2}\Big(-\a s^n + 2\Big); $$ 
\m
\item {(d)} $A = A(2,-1,0)$ and $\l_n = 0$ for all $n$;
\m
\item {(e)} $A = A(0,0,\ga)$ and $\l_n = \ga$ for all $n \geq 1$.   
\endproclaim
\b
{\bf Proof.} We will divide considerations into cases according
to the solutions given for $\l_n$  by Proposition 2.12.
\m
\n {\it Case (a):  $\a+4 \be \neq 0$ and $\a + \be \neq 1$.}
\m
In this first case $r_1 \neq r_2$ and 

$$\l_n = c_1 r_1^n + c_2 r_2^n + \frac {\ga}{1 - \a -\be}.$$

\n From $\l_\ell = \l_k$ and $\l_{\ell-1} = \l_{k-1}$ we see

$$\aligned & c_1 r_1^k (r_1^{\ell-k} -1) = - c_2 r_2^k(r_2^{\ell-k} -1) \\
 & c_1 r_1^{k-1} (r_1^{\ell-k} -1) = - c_2 r_2^{k-1}(r_2^{\ell-k} -1). \\
\endaligned$$

\n Multiplying the second by $r_2$ and subtracting gives

$$c_1r_1^{k-1}(r_1^{\ell-k}-1)(r_1-r_2) = 0.$$

\n This shows that $c_1 = 0, r_1 = 0$, or $r_1^{\ell-k} = 1$.
Analogously, $c_2 = 0, r_2 = 0$, or $r_2^{\ell-k} = 1$.  

Suppose initially that $0 = r_2 = (1/2)\a - t$.  Then $r_1 = (1/2)\a + t = \a
\neq 0$, and $t^2 = (1/4)\a^2$ implies $\be = 0$.   From Proposition 2.12 we
see

$$\l_n = c_1 \a^n + \frac {\ga}{1-\a} = \Big(\l - \frac
{\ga}{(1-\a)}\Big) \a^n + \frac {\ga}{1-\a}. \tag 2.14$$

\n  When $\ga = 0$,
then $\l_n = \l \a^n$, and $\l_\ell = \l_k$ forces $\l = 0$ (Case (1))
or $\a^{\ell-k} = 1$ (Case (3) with $\ga = 0$).   When $\ga \neq 0$, then 
either $\l_n = (1-\a)^{-1}\ga$ for all $n$ (Case (2))
or $\l_n$ is as in (2.14) where $\a^{\l-k} = 1$ (Case (3) with $\be = 0$).  

When $r_1 = 0$ and $r_2 = \a$, the solutions are 
precisely the same.  So we can assume that neither $r_1$ nor $r_2$
is 0.  

Suppose next that $c_1 = 0$.  Then $\l = \l_0 = c_2 + (1-\a-\be)^{-1}\ga$,  
which gives $\l_n =
\Big(\l-(1-\a-\be)^{-1}\ga\Big)r_2^n + (1-\a-\be)^{-1}\ga$. The condition
$c_1 = 0$ says that $r_1 \l = (1-\a-\be)^{-1}\ga(\a+\be-r_2)$.
When  $c_2 =0$, then $\l = (1-\a-\be)^{-1}\ga$ and we have 

$$\displaystyle\frac {r_1 \ga}{1-\a-\be}= 
\displaystyle\frac {(\a+\be-r_2) \ga}{1-\a-\be}.$$

\n If $\ga = 0$, then $\l = 0$ and $\l_n = 0$ for all $n$ (Case (1)).
If $\ga \neq 0$, then $\a = r_1+r_2 = \a+\be$ forces $\be = 0$
and $\l_n = (1-\a)^{-1}\ga$ (Case (2)).  The only other option is
that $c_2 \neq 0$ but $r_2^{\ell-k} = 1$.  That gives Case (3).
The situation when $c_2 = 0$ and $c_1 \neq 0$ also gives Case (3).  

Finally, if $r_1,r_2,c_1,c_2$ are all nonzero, then
$r_1^{\ell - k} = 1 = r_2^{\ell-k}$ and

$$\l_n = c_1 r_1^n + c_2 r_2^n + \displaystyle\frac {\ga}{1-\a-\be},$$

\n for all $n$ (Case (4)).  

\m

{\it Case (b):  $\a+4\be \neq 0$ and $\a+\be = 1$.}
\m
Under these hypotheses, $\be = 1-\a$ and 
$0 \neq \a^2 + 4(1-\a) =(2-\a)^2$.  Then
$t^2 = (1/4)(2-\a)^2$, and we may assume that $t = -(1/2)(2-\a)$.  From that 
we have $r_1 = \a -1$, $r_2 = 1$, and

$$\l_n = c_1(\a-1)^n + c_2 + \displaystyle \frac {\ga n}{2-\a}. \tag 2.15$$

\n When $\l_\ell = \l_k$ and $\l_{\ell-1} = \l_{k-1}$, then

$$\aligned 
& c_1(\a-1)^\ell + c_2 + \displaystyle \frac {\ga \ell}{2-\a}
= c_1(\a-1)^k + c_2 + \displaystyle \frac {\ga k}{2-\a} \\
& c_1(\a-1)^{\ell-1} + c_2 + \displaystyle \frac {\ga (\ell-1)}{2-\a}
= c_1(\a-1)^{k-1} + c_2 + \displaystyle \frac {\ga (k-1)}{2-\a}\\
\endaligned$$

\n or

$$
\aligned
c_1(\a-1)^k \Big((\a-1)^{\ell-k} -1\Big) =  \displaystyle \frac {\ga
(k-\ell)}{2-\a}
\\
c_1(\a-1)^{k-1} \Big((\a-1)^{\ell-k} -1\Big) =  \displaystyle \frac {\ga
(k-\ell)}{2-\a}.
\\
\endaligned \tag 2.16$$

\n Therefore,  

$$c_1(\a-1)^{k-1}(\a-2)\Big((\a-1)^{\ell-k} -1\Big) = 0.$$

\n Since $\a \neq 2$, we see that the left-hand sides of (2.16) are 0, and
$\ga = 0$.  Solving for $c_1,c_2$ from Proposition 2.12 we find that $c_1 =
(2-\a)^{-1}(1-\a)\l$ and $c_2 = (2-\a)^{-1}\l$.   Therefore,

$$\l_n = \frac {\l(1-\a)}{2-\a}\,(\a-1)^n + \frac {\l}{2-\a}\,1^n = 
\frac {\l}{2-\a}\Big( 1 - (\a-1)^{n+1}\Big). \tag 2.17$$

\n The condition $\l_{\ell} = \l_k$ implies that $\l = 0$, and hence
all $\l_n = 0$ (Case (1)), or that $\a = 1$ 
and hence that $\l_n = \l$ for all $n$ (Case (1)), or that $\l \neq 0$ and
$(\a-1)^{\ell-k} = 1$ (Case (2)).   
\b
\n {\it Case (c):  $\a^2 + 4 \be = 0, \; \a + \be \neq 1$ and $\a \neq 0$.}
\m
Here $\be = -(1/4)\a^2$ so that

$$\l_n = c_1 s^n + c_2 n s^n + \displaystyle \frac {\ga}{1-\a-\be} 
= c_1 s^n + c_2 n s^n + \displaystyle \frac {4\ga}{(2-\a)^2}.
\tag 2.18$$ 

\n It follows from $\l_\ell= \l_k$ and $\l_{\ell-1} = \l_{k-1}$ that

$$
\aligned
& c_1s^k\Big(s^{\ell-k}-1\Big) = -c_2 s^k \Big (\ell s^{\ell-k} - k\Big) \\
& c_1s^{k-1}\Big(s^{\ell-k}-1\Big) = -c_2 s^{k-1} \Big ((\ell-1) s^{\ell-k} - 
(k-1)\Big) \\
\endaligned \tag 2.19$$

\n Since $s = (1/2)\a \neq 0$, we may use these equations to see that 

$$c_2(s^{\ell-k} - 1) = 0 = c_1 (s^{\ell-k}-1).$$

\n Now if $s^{\ell-k} = 1$, then (2.19) shows that $c_2 = 0$, which
says that $\l  = (2-\a)^{-1}(2\ga)$.  Therefore, since $c_1 =
\l-(2-\a)^{-2}4\ga$, we have

$$\l_n = \Big(\l - \frac {4\ga}{(2-\a)^2}\Big)s^n + \frac {4\ga}{(2-\a)^2} =
\frac {2\ga}{(2-\a)^2}\Big(-\a s^n + 2\Big), \tag 2.20$$

\n where $s = (1/2)\a$ satisfies $s^{\ell-k} = 1$.  This gives the solutions in (c).
Now if $s^{\ell-k} \neq 1$, then $c_1 = 0 = c_2$ and
$\l_n = (2-\a)^{-2}4\ga$ for all $n$.   Those relations imply
$(2-\a)^{-1}2\ga = (2-\a)^{-2}4\ga$, which gives $\ga = 0$
or $\a = 0$.  Since $\a \neq 0$ by assumption, $\ga = 0$ and $\l_n = 0$
for all $n$.    
 
\b
\n {\it Case (d):  $\a^2 + 4 \be = 0$ and  $\a = 2$ and $\be = -1$.} 
\m
The values of $\l_n$ are given by
$\l_n = c_1 s^n + c_2 n s^n + (1/2)\ga n^2$, which reduces to
$\l_n = c_1 + c_2 n + (1/2)\ga n^2$  because $s = \a/2 = 1$.  Solving
for $c_1$ and $c_2$ using Proposition 2.12, 
 we see that $c_1 = \l$ and $c_2 = \l + \ga/2$. 
Therefore, 

$$\aligned
\l_n & = c_1 + c_2 n + \displaystyle \frac {\ga n^2}{2} 
 = \l + (\l+\displaystyle \frac {\ga}{2})n  + \displaystyle \frac {\ga n^2}{2}
\\  & = (\displaystyle \frac {\ga}{2}n + \l)(n+1). \\
\endaligned \tag 2.21$$

\n Then
 
$$\l_{\ell} = \l + (\l+\displaystyle \frac {\ga}{2})\ell  + \displaystyle \frac {\ga
\ell^2}{2}
= \l + (\l+\displaystyle \frac {\ga}{2})k  + \displaystyle \frac {\ga k^2}{2}
= \l_k$$

\n implies that  

$$(\displaystyle \frac {\ga}{2})(\ell^2 - k^2) + (\l+\displaystyle \frac
{\ga}{2})(\ell-k) = 0.$$

\n Similarly, from $\l_{\ell-1} = \l_{k-1}$ we obtain 

$$(\displaystyle \frac {\ga}{2})((\ell-1)^2 - (k-1)^2) + (\l+\displaystyle \frac
{\ga}{2})(\ell-k) = 0,$$

\n and those equations combine to give

$$\ga(\ell-k) = 0.$$

\n It must be $\ga = 0$ and $\l_n = \l(n+1)$.  But then 
$\l_\ell = \l(\ell+1) = \l(k+1) = \l_k$ shows that $\l = 0$.
We conclude that for the algebra $A(2,-1,\ga)$, the only
time two weights $(\l_\ell,\l_{\ell-1})$ and  $(\l_k,\l_{k-1})$
can be equal is if $\ga = 0$ and $\l_n = 0$ for all $n$. 

\m
\n {\it Case (e):  $\a = \be = 0$}
\m
In this case $\l_n = \ga$ for all $n \geq 1$.    \qed
\b
\proclaim{Corollary 2.22} In $V(\l)$ each weight space is either
one-dimensional or infinite-dimensional.  If an infinite-dimensional
weight space occurs, there are only finitely many weights.  The cases when 
$V(\l)$ has an
infinite-dimensional weight space are the ones listed in
Theorem 2.13.  \endproclaim   
\b
\m
\subhead Submodules of $V(\l)$ \endsubhead
\b
\proclaim{Proposition 2.23} If each weight space of $V(\l)$ is
one-dimensional (i.e. if we are not in one of the cases described
by Theorem 2.13), then the proper submodules of $V(\l)$ have the
form $N =$ span$_\C\{v_j \mid j \geq n\}$ for some $n > 0$ with $\l_{n-1} = 0$,
and hence they are contained in $M(\l)$.  \endproclaim
\b
{\bf Proof.} Suppose $N$ is a proper submodule
of $V(\l)$.  There exists some weight $\nu = (\l_n,\l_{n-1})$ such
that $N_\nu \neq (0)$, and we may suppose that $n$ is the
least value with that property.   Then $v_n \in N$, and $n$ is
the smallest value for which that is true.   Applying $u$ we
see that $v_j \in N$ for all $j \geq n$.   \qed  
\b
To describe the submodules of $V(\l)$ when
its weights are given by Theorem 2.13,  
it is helpful to observe that if
$N$ is a proper submodule of any $V(\l)$, then

$$I_N = \{f \in \C[x] \mid f(u) \cdot v_0 \in N\} \tag 2.24$$

\n is an ideal of $\C[x]$.   There exists a monic polynomial
$g(x)$ which generates that ideal.   Then 
$\dim_\C V(\l)/N
= \deg(g) = k$, and $v_i+N$ for $i = 0,1, \dots, k-1$ is a basis
for $V(\l)/N$.
\b
The situations described by Theorem 2.13 can be grouped into
four general types (see 2.17 and 2.20): 
\b
\n (2.25)
 \item {{}(i)} $\ga = 0$ and $\l_n = 0$ for all $n$.  

\item{{}(ii)} $\l_n = c$ for all $n$ where $c \in \C$ satisfies $c-\a c - \be c = \ga$.  

\item {{}(iii)}  $\l_n = (\l-c)r^n + c$  where $c \in \C$ and $r$ is a root
of $x^2 - \a x - \be = 0$.  There is a least integer $p \in \Z_{>0}$  such that $r^p = 1$.  

\item {{}(iv)}  $\l_n = c_1 r_1^n + c_2 r_2^n + c$ where $r_1 \neq r_2$ are
roots of the equation $x^2 - \a x - \be = 0$, and there is some smallest
value $p \in \Z_{>0}$ so that $r_1^p = 1 = r_2^p$.  
\b
(i) Assume first that $A = A(\a,\be,0)$ and consider the module $V(0)$. 
Since $\l_n = 0$  for all $n$, for any $g \in \C[x]$, 

$$N^{(g)} = \text{span}_\C\{f(u)g(u)v_0 \mid f \in \C[x]\} \tag 2.26$$ 

\n is a submodule of $V(0)$, which is proper if $\deg(g) \geq 1$.  
For each $\xi \in \C$, the submodule $N^{(x-\xi)}$ is maximal, and
every maximal submodule has this form.  In particular, $M(\l) = N^{(x)}$.
The quotient module $L(0,\xi) \eqdef V(0)/N^{(x-\xi)}$ is spanned by $v_0 + N^{(x-\xi)}$
where $d \cdot (v_0 + N^{(x-\xi)}) = 0$ and $u \cdot (v_0 + N^{(x-\xi)}) = 
\xi v_0 + N^{(x-\xi)}.$
\b
(ii) Next suppose that $A = A(\a,\be,\ga)$, and the module is $V(c)$ where
$0 \neq c \in \C$ is such that $c -\a c -\be c = \ga$.  Then
$\l_n = c \neq 0$ for all $n$.  By Proposition 2.4, $V(c)$ is simple. 
\b
(iii) and (iv) These two cases behave similarly, and so
we will treat them simultaneously.  Assume for $A = A(\a,\be,\ga)$ that the
module $V(\l)$ has weights given by   

$$
\aligned
& \l_n = (\l-c)r^n + c \; \text {where} \;  c \in \C, \; \l \neq c, \; \text {and} \; r^p = 1, \;
\text {or} \\
& \l_n  = c_1 r_1^n + c_2 r_2^n + c \; \text {where} \; c \in \C, \;
c_1, c_2 \neq 0, \; \text {and}
\; r_1^p = 1 = r_2^p.
\\
\endaligned $$

\n In each instance, we will assume that $p$ is chosen to be
minimal.    Then $\l_{n+p} = \l_n$ for all $n \geq 0$, and $V(\l) = \sum_{i = 0}^p
V(\l)_{\nu_i}$, where
$\nu_i = (\l_i,\l_{i-1})$.  Now $\nu_p = \nu_0$ if and only if $\l_{p-1} = 0$, which happens, for
example, if
$\be \neq 0$. 

Suppose first that $\nu_p \neq \nu_0$, so that $V(\l)$ has $p+1$ distinct weights.
If $N$ is a submodule of $V(\l)$ and $N_{\nu_0} \neq (0)$, then $N = V(\l)$. 
In particular, any proper submodule $N$ is contained in $M(\l)$ by Proposition 2.4,
so that $M(\l)$ is the unique maximal submodule. 
We may assume that $N_{\nu_j} \neq (0)$ for some $j = 1, \dots, p$.  Then
$(0) \neq u^{p-j+k}N_{\nu_j} \subseteq N_{\nu_k}$ so that all $N_{\nu_k} \neq (0)$.
Therefore,  for $k \neq p$, 

$$\{f \in \C[x] \mid u^k f(u^p) \cdot v_0 \in N_{\nu_k}\} \tag 2.27$$

\n is a nonzero ideal of $\C[x]$, and we may assume that $g_k$ is its monic generator.
Similarly, $\{f \in \C[x] \mid f(u^p) \cdot v_0 \in N_{\nu_p}\}$ is a nonzero
ideal, whose monic generator we denote by $g_p$.  Then $u^k g_k(u^p) \cdot v_0 \in
N_{\nu_k}$ implies $u (u^k g_k(u^p) \cdot v_0) \in N_{\nu_{k+1}}$.   It follows that

$$g_{p-1} \mid g_{p-2} \mid \dots \mid g_2 \mid g_1 \mid g_p \mid xg_{p-1}.$$

\n Either $g_1 = g_2 = \dots = g_p = g$ and $N = N^{(g(x^p))}$ 
(see (2.26) for the notation), or
there is some value $m$ such that $g_{p-1} = \dots = g_{m+1} = g$ and
$g_m = \dots = g_1 = g_p = xg$.  Since $d (u^{m+1}g(u^p) \cdot v_0)
= \l_m u^m g(u^p) \cdot v_0,$ it must be that $\l_m = 0$ in this case.
It is easy to see that $\l_k \neq \l_\ell$ for $1 \leq k < \ell \leq p$.  
Therefore $m$ is the minimal value such that $\l_m = 0$. In this
situation, $N = N^{(x^{m+1}g(x^p))}$.  As a special case, we
have $M(\l) = N^{(x^{m+1})}$.

In the case that  $\nu_p = \nu_0$, the module $V(\l)$ has $p$ distinct weights.  
Assume $N$ is a proper submodule of $V(\l)$. 
Arguing as above,  we have that $N_{\nu_k} \neq (0)$ for each $k = 0,1, \dots, p-1$.
Letting $g_k$ be the monic generator for the ideal in (2.27), we have

$$g_{p-1} \mid g_{p-2} \mid \dots \mid g_2 \mid g_1 \mid g_0 \mid xg_{p-1}.$$

\n Again there are two possibilities.  Either $g_0 = g_1 = \dots = g_{p-1} = g$
and $N = N^{(g(x^p))}$, or else 
there is some value $m$ such that $g_{p-1} = \dots = g_{m+1} = g$ and
$g_m = \dots = g_1 = g_0 = xg$.   Here
$N = N^{(x^{m+1}g(x^p))}$, and  $M(\l) = N^{(x^{m+1})}$.
\b
\proclaim {Corollary 2.28} If $\gamma = 0 = \lambda$, then
$V(\l)$ has infinitely many maximal proper submodules, each
of the form 
$N^{(x-\xi)} = \text{span}_\C\{f(u)(u-\xi)v_0 \mid f \in \C[x]\}
= \text{span}_\C\{v_n -\xi v_{n-1} \mid n = 1,2,\dots\}$
for some $\xi \in \C$, and infinitely many one-dimensional
simple modules $L(0,\xi) = V(0)/N^{(x-\xi)}$.  In
all other cases, $M(\l)$ is the unique maximal submodule of $V(\l)$,
and there is a unique simple highest weight module
of weight $\l$, $L(\l) \eqdef V(\l)/M(\l)$, up to
isomorphism.   \endproclaim   
 
\b
\subhead Enveloping algebra examples \endsubhead
\b
Recall that the universal enveloping algebra $U(\fsl)$ of $\fsl$ is
isomorphic to the algebra $A(2,-1,-2)$, and the universal enveloping
algebra $U({\frak H})$ of the Heisenberg Lie algebra $\frak H$ is isomorphic
to $A(2,-1,0)$.   For all the 
algebras
$A(2,-1,\ga)$, we have computed the values of $\l_n$ in 
the proof of Case (d) of Theorem 2.13, and from (2.21) we see that

$$\l_n = (\displaystyle \frac {\ga}{2}n + \l)(n+1).$$
 
In the $\fsl$-case, the operator $h = du-ud$ is used rather
than $du$.  The eigenvalues of $h$ are $\l_n -\l_{n-1} = \l+n\ga = \l-2n$,
$n = 0,1, \dots$, 
(as is customary in the representation theory of $\fsl$), and
$V(\l)$ is simple if and only if $\l \not \in \Z_{\geq 0}$.   The analogous
computation in the Heisenberg Lie algebra shows that the central element $z = du-ud$
has  constant eigenvalue $\l_n = \l$.  
\b
\m
\subhead Quantum examples \endsubhead
\b
Let
$\g$ be a symmetrizable Kac-Moody Lie algebra (for example, a 
finite-dimensional 
simple Lie algebra) over $\C$ corresponding
to the Cartan matrix  $\frak A = (a_{i,j})_{i,j= 1}
^n$.  Then there are relatively prime integers $\ell_i$ so that
the matrix $(\ell_i a_{i,j})$ is symmetric. 
Let $\C(q)$ be the field of rational functions in the
indeterminate $q$ over $\C$.  Assume

$$q_i = q^{\ell_i},\quad \quad \text {and} \quad \quad  
[m]_i = \frac {q_i^{m} - q_i^{-m}} {q_i -
q_i^{-1}}$$

\n for all $m \in \Z_{\geq 0}$.  When $m \geq 1$,  let 
$$[m]_i! = \prod_{j = 1}^m [j]_i.$$

\n Set $[0]_i! = 1$ and define

$$\left [\matrix m \\ n \\ \endmatrix
\right]_i = \frac{[m]_i!}
{[n]_i! [m-n]_i!}.$$ 
\n Then the quantized enveloping algebra (quantum group) $U_q(\g)$
of $\g$ is the unital associative algebra over $\C(q)$ with generators 
$E_i, F_i, K_{i}, K_{i}^{-1}$ ($i = 1, \dots, n$) subject to the relations
\b
\item {}{(Q1)} $K_i K_i^{-1} = K_i^{-1} K_i,$
\quad \quad $K_i K_j = K_jK_i$
\m
\item {}{(Q2)} $K_i E_j K_i^{-1} =
q_i^{a_{i,j}}E_j$ 
\m
\item {}{(Q3)} $K_i F_j K_i^{-1} =
q_i^{-a_{i,j}}F_j$ 
\m
\item {}{(Q4)} $E_i F_j - F_j E_i
= \displaystyle \delta_{i,j} \frac {K_i - K_i^{-1}}
{q_i - q_i^{-1}}$ 
\m
\item {}{(Q5)}$\displaystyle \sum_{k = 0}^{1-a_{i,j}} (-1)^k
\left [\matrix 1-a_{i,j}\\ k  \\ \endmatrix
\right]_i E_i^{1-a_{i,j}-k} E_j E_i^k = 0$ \quad \quad  for  \quad $i \neq j$
\m
\item {}{(Q6)}$\displaystyle \sum_{k = 0}^{1-a_{i,j}}
(-1)^k\left [\matrix 1-a_{i,j}\\ k  \\ \endmatrix
\right]_i F_i^{1-a_{i,j}-k} F_j F_i^k = 0$ \quad \quad for \quad $i \neq j.$
\b
Suppose $a_{i,j} = -1 = a_{j,i}$ for some $i \neq j$, and consider the subalgebra
$U_{i,j}$ generated by $E_i, E_j$. In this special case, the quantum Serre
relation (Q5) reduces to

$$\aligned & E_i^2 E_j - [2]_i E_i E_j E_i +
E_j E_i^2 = 0\quad \quad \text {and}  \\
& E_j^2 E_i - [2]_j E_j E_i E_j +
E_i E_j^2 = 0.\\
\endaligned $$

\n Since $-\ell_i = \ell_ia_{i,j} = \ell_ja_{j,i} = -\ell_j$, the coefficients
$[2]_i$ and $[2]_j$ are equal. 
The algebra  $U_{i,j}$ (with $q$ is specialized to
a complex number which is not a root of unity) is isomorphic to
$A([2]_i,-1,0)$ by the mapping $E_i \mapsto d$,
$E_j \mapsto u$.   The same result is true if
the corresponding $F$'s are used in place of the $E$'s. 
In particular, when $\g = \fst$, the 
algebra $U_{i,j}$  is just the subalgebra of $U_q(\fst)$
generated by the $E$'s. 
 
\b
Here we compute the values of $\l_n$ in the
case of the algebra $A([2]_i,-1,0)$. For convenience
of notation we write 

$$p = q_i.$$ 

\n Then
$\a = [2]_i = \displaystyle \frac {p^2 - p^{-2}} {p - p^{-1}} = p+
p^{-1}$,
$\be = -1$ and $\ga = 0$ so that

$$\a^2  + 4 \be = p^2 + 2 + p^{-2} -4 = (p - p^{-1})^2.$$

\n Thus 

$$\aligned r_1 & = \frac {p + p^{-1} +  p-p^{-1}} {2} = p \\
r_2 & = \frac {p + p^{-1} - (p-p^{-1})} {2} = p^{-1}, \\
\endaligned$$

$$\aligned \left (\matrix c_1 \\ c_2 \\
\endmatrix
\right) & = \frac {1}{p^{-1}-p}
\left (\matrix p^{-1} & -1 \\
-p & 1 \\ \endmatrix
\right) \left (\matrix \l \\
(p+p^{-1})\l \endmatrix
\right) \\ &  = \frac {1}{p^{-1}-p}
\left (\matrix -p\l \\
p^{-1} \l \endmatrix
\right) = 
\frac {1}{p-p^{-1}}
\left (\matrix p \l \\ -p^{-1}\l
\endmatrix
\right). \\ \endaligned $$

\n Therefore 

$$\aligned \l_n & = \frac {p\l}{p-p^{-1}}p^n - \frac
{p^{-1}\l}{p-p^{-1}}p^{-n} \\
& = \Big(\frac {p^{n+1} - p^{-(n+1)}} {p-p^{-1}}\Big)\l \\
& = [n+1]_i \l. \\
\endaligned
$$
\b
\subhead $(q,r)$-differential poset examples 
\endsubhead
\b For an algebra arising from a $(q,r)$-differential poset 
the case of interest is the
algebra $A= A(\a,\be,\ga)$ with defining relations

$$\aligned
d^2 u & = q(q+1)dud -q^3ud^2 + rd, \\
du^2 & = q(q+1)udu -q^3u^2d + ru. \\
\endaligned $$

\n Here $r$ is an arbitrary element of the field.
For this algebra $\a^2 + 4\be
= q^2(q-1)^2$, so we may assume that $t = (1/2)q(q-1)$ and 
$r_1 = 1/2(\a + q(q-1)) = q^2$ and 
$r_2 = 1/2(\a - q(q-1)) = q$.  
(We will suppose throughout that $q \neq 0, \pm 1$).  
Then 

$$1-\a -\be = 1-q(q+1)+q^3 = q^2(q-1)-(q-1) = (q-1)(q^2-1)$$

\n so that

$$x_n = (1-\a-\be)^{-1}\ga = \frac {r}{(q-1)(q^2 -1)}$$

\n for all $n$.  We call that common value $x$. Recall that

$$\left (\matrix c_1 \\ c_2 \\ \endmatrix \right )
 = \frac {1}{r_2-r_1}
\left (\matrix r_2 & -1 \\ -r_1 & 1 \\ \endmatrix \right )
\left (\matrix \l-x \\ \a\l+r-x \\ \endmatrix \right ). $$

\n Now 

$$\frac {1}{r_2-r_1} = \frac {1}{q(1-q)},$$

\n and putting that into the above equation, we get after
a bit of calculation

$$\aligned c_1 & = \frac {1}{q(1-q)}\Big(q(\l-x) - (\a\l+r -x)\Big) \\
& = \frac {-q}{1-q}\Big(\l+\frac{r}{q^2-1}\Big) \\
c_2 & = \frac {1}{q(1-q)}\Big(-q^2(\l-x) + (\a\l+r -x)\Big) \\
& = \frac {1}{1-q}\Big(\l+\frac{r}{q-1}\Big). \\ 
\endaligned $$

\n Thus the eigenvalues of $du$ for this algebra are

$$\aligned
\l_n & = \frac {-q}{1-q}\Big(\l+\frac{r}{q^2-1}\Big)q^{2n}
+\frac {1}{1-q}\Big(\l+\frac{r}{q-1}\Big)q^n
+ \frac {r}{(q-1)(q^2 -1)} \\
& = \frac{q^{n+1}-1}{q-1}\Big( q^n \l + \frac {q^n-1}{q^2-1} r\Big).
\endaligned $$ 
\b
\subhead Fibonacci examples 
\endsubhead
\b For the algebra $A(1,1,0)$, 
 the solutions to the associated
linear recurrence $\l_n = \l_{n-1} + \l_{n-2}$, $\l_0 = \l$, $\l_{-1} = 0$,
(hence the eigenvalues of $du$ on $V(\l)$)   
are given by the Fibonacci sequence
$\l_0 = \l$, $\l_1 = \l$, $\l_2 = 2 \l$, $\l_3 = 3\l$, $\l_4 = 5 \l$,
$\dots$.  In this case, the equations in Proposition 2.12 reduce to

$$\l_n = c_1\left (\frac {1 + \sqrt 5} {2}\right)^n + c_2 \left (\frac {1 - \sqrt 5}
{2}\right)^n,$$

\n where 

$$ \aligned 
c_1 & = \frac {\l \sqrt 5}{5}\left (\frac {1 + \sqrt 5} {2}\right ) \\
c_2 & = -\frac {\l \sqrt 5}{5}\left(\frac {1 - \sqrt 5} {2}\right) \\
\endaligned $$

\n or  

$$\l_n = \l\frac {\sqrt 5}{5}\left( \left (\frac {1 + \sqrt 5} {2}\right)^{n+1}
-\left(\frac {1 - \sqrt 5} {2}\right)^{n+1}\right).$$

\b
\subhead  Lowest weight modules \endsubhead 
\b 
To create lowest weight modules for $A = A(\a,\be,\ga)$
we reverse the roles of $d$ and $u$.  Thus, we assume there
is a vector $w_0$ such that $u \cdot w_0 = 0$, 
$ud \cdot w_0 = \k w_0$ and $W = Aw_0$.    
When $\be \neq 0$, the  eigenvalues of $du$ are given by
the sequence which has $\k_{-1} = 0$, $\k_0 = \kappa$,
an arbitrary complex number, and 

$$\aligned & \be \k_n + \a \k_{n-1} - \k_{n-2} = -\ga, \quad \text {or equivalently}
\\ & \k_{n} = \be^{-1}\big(-\a \k_{n-1} + \k_{n-2} -  \ga\big)
 \\
\endaligned 
\tag 2.29$$ 
\n for all $n \geq 1$.  
\b
\proclaim {Proposition 2.30} 
 Let $W(\k)$ be the
$\C$-vector space having basis $\{w_n \mid n = 0,1,2, \dots\}$.
\s
\item {(a)} Assume $\be \neq 0$,  $\k_{-1} = 0$, and  
$\k_0 = \kappa$, an arbitrary element of $\C$.  Suppose
$\k_n$ for $n \geq 1$ is as in (2.29), and define

$$\aligned
u \cdot w_n & = \k_{n-1} w_{n-1}, \quad n \geq 1, \quad \text {and}
\quad u
\cdot w_0 = 0 \\
d \cdot w_n & = w_{n+1}. \\
\endaligned \tag 2.31$$

\item {}Then this action gives $W(\k)$ the structure of
a lowest weight $A(\a,\be,\gamma)$-module.
\m
\item {(b)} When $\be = 0$ and $\a \neq 0$, let $\k_n$
be
defined by 

$$\k_n = -\ga \sum_{j = 1}^{n+1} \a^{-j}$$

\item {}for all $n \geq 0$.  Then  $W(-\ga \a^{-1})$ with the
action given by (2.31) is a lowest weight $A(\a,\be,\ga)$-module. 
\m
\item {(c)} When  $\ga \neq 0$, the only lowest
weight $A(0,0,\ga)$-module is the 1-dimensional module $W = \C w_0$ 
with $d \cdot w_0 = 0 = u \cdot w_0$.  When $\ga = 0$, set $\k_n = 0$
for all
$n$. Then
$W(0)$ with the action given by (2.31) is a lowest weight $A(0,0,0)$-module.     
\endproclaim
\b 
The proof  of 
this result is identical to that of Proposition 2.2, and the following
proposition can be shown using standard recurrence relation techniques.   
 
\b

\proclaim { Proposition 2.32}
Assume $\be \neq 0$ and fix  $t \in \C$ such that

$$t^2 = \frac {\a^2 + 4 \be}{4\be}.$$

\n Let  $\k_{-1} = 0$, $\k_0 = \kappa \in \C$,   
and for
$n \geq 1$, let $\k_n$  be given by the recurrence relation
$\be \k_n+\a\k_{n-1} -\k_{n-2} = -\ga$.
\m
\item {}{(i)} If  $\a^2 + 4 \be\neq 0$, then 

$$\k_n = c_1 r_1^n
+ c_2 r_2^n  + x_n,$$

\n where 

$$\aligned & r_1 = -\frac {\a} {2\be} + t, \quad
 \quad 
 r_2 = -\frac {\a}{2\be} -t , \\
& x_n = \cases 
(1-\a-\be)^{-1}\ga \ \ \text {if $\a+\be \neq 1$} \\
(\a-2)^{-1}\ga n \ \ \text {if $\a+\be = 1$ and $\a \neq 2$,} \\
\endcases \\
 & \text {and} \quad \quad \left (\matrix c_1 \\ c_2 \\ \endmatrix \right )
 = \frac {1}{r_2-r_1}
\left (\matrix r_2 & -1 \\ -r_1 & 1 \\ \endmatrix \right )
\left (\matrix \l-x_0 \\ \a\l+\ga-x_1 \\ \endmatrix \right ). 
\endaligned $$
 \m
\item {}{(ii)} If $\a^2 + 4 \be = 0$, then $\a \neq 0$ and

$$\k_n = c_1 s^n + c_2 n s^n + x_n $$

\n where  

$$\aligned 
& s = \displaystyle \frac {\a}{2\be} \\
& x_n = \cases 
(1-\a-\be)^{-1}\ga \ \ \text {if $\a+\be \neq 1$} \\
\displaystyle -\frac {\ga}{2} n^2 \ \ \text {if $\a+\be = 1$ i.e. if  $\a \neq 2$ and $\be = -1$,} \\
\endcases \\
& {and}  \quad \quad \left (\matrix c_1 \\ c_2 \\ \endmatrix \right)
 =  
\left (\matrix 1 & 0 \\ \displaystyle -\frac{2\be}{\a} & \displaystyle \frac {2\be}{\a} \\ \endmatrix
\right )
\left (\matrix \k-x_0 \\ \displaystyle -\frac{(\ga + \a \k)}{\be}-x_1 \\ \endmatrix
\right ). 
\endaligned $$

\endproclaim
\b 
\subhead Doubly-infinite modules \endsubhead
\b 
\proclaim {Proposition 2.33} Assume $\be \neq 0$ and let  $\k,\l$
be arbitrary but fixed complex numbers.  Suppose $\l_0 = \l$, $\l_1 = \a \l + \be \k
+\ga$, and for 
$n \geq 2$, $\l_n$ is defined by the recurrence relation $\l_n =
\a \l_{n-1} + \be \l_{n-2} + \ga$.  Set $\l_{-1} = \k$ and define $\l_{n-2}$
for all $n \leq 0$ by the equation 
$\l_{n} = \a \l_{n-1}+ \be \l_{n-2}+\ga$.
The complex vector space $V(\k,\l)$
with basis $\{v_n \mid n \in \Z\}$ and with

$$\aligned & d \cdot v_n = \l_{n-1} v_{n-1}  \\
& u \cdot v_n =  v_{n+1} \\
\endaligned  \tag 2.34$$

\n is an $A(\a,\be,\ga)$-module.  \endproclaim 
\b
{\bf Proof.}  As in our previous arguments, we may assume there is
an action of the free associative algebra generated by $d$ and $u$ 
on $V(\k,\l)$ given by (2.34).  Then 

$$\aligned \Big (d^2 u - \a dud - \be u d^2 -\ga d \Big)v_m & =
\l_{m-1}\Big(\l_m- \a \l_{m-1}- \be\l_{m-2} -\ga \Big)
v_{m-1} = 0
\\
\Big (d u^2 - \a udu - \be u^2 d - \ga u \Big)v_m
& = \Big(\l_{m+1} - \a \l_m - \be \l_{m-1} - \ga\Big)v_m = 0, \\
\endaligned $$

\n so that there is an induced action by $A(\a,\be,\ga)$.  \qed 
\b
If $\l_m = 0$ for any $m$, then the span of the vectors $\{v_n \mid n \geq m+1\}$
is a submodule of $V(\k,\l)$. 
\b 
\subhead Dual modules \endsubhead
\b Suppose $V$ is any module for $A = A(\a,\be,\gamma)$ and let
$\{v_a \mid a \in I\}$ be a basis for $V$.  Assume $v_a^* \in \Hom_\C(V,\C)$
is defined by $v_a^*(v_b) = \delta_{a,b}$ for all $a,b \in I$.
Let $V^*$ be the finite dual of $V$ (i.e. the finite linear combinations
of the elements $v_a^*$).   We will show that $V^*$ carries the
structure of an $A$-module. 
\b
\n {\it The antiautomorphism $\eta$}
\m
Define a mapping $\eta: A \rightarrow A$ by first
specifying
$\eta(u) = d$ and $\eta(d) = u$ and extending this to
an antiautomorphism of the free associative
algebra generated by $u$ and $d$.  Then

$$\eta(d^2u) = du^2 = \a udu + \be u^2 d + \ga u 
= \eta(\a dud + \be u d^2 + \ga d),$$

\n and a similar relation holds for $du^2$.  Therefore $\eta$
induces an antiautomorphism of order 2 on $A$ (which we also
denote by $\eta$).  
\b
Now define 

$$(x \cdot v^*)(w) = v^*(\eta(x)\cdot w)$$

\n for all $x \in A$, $w \in V$, and $v^* \in V^*$.  It is easy to see that
this makes $V^*$ an $A$-module. 
\b
Suppose now that $V = V(\l)$ is a highest weight $A$-module
with basis $\{v_n \mid n = 0,1, \dots\}$.  Then

$$\aligned
(u \cdot v_\ell^*)(v_n) & = v_\ell^*(d \cdot v_n) = \l_{n-1} v_\ell^*(v_{n-1}) \\
& = \cases  
0 \ \ \text {if $n \neq \ell+1$} \\ 
\lambda_\ell \ \ \text {if $n = \ell+1.$} \\
\endcases
\\ 
\endaligned$$

\n Therefore $u \cdot v_\ell^* = \l_\ell v_{\ell+1}^*$.  Similarly,
$d \cdot v_\ell^* = v_{\ell-1}^*$.   If $\l_n = 0$ for some $n$, then
$V(\l)^*$ is not generated by $v_0^*$.   But whenever $\l_n \neq 0$ for all $n$,
then $V(\l)^*$ is a highest weight module.  In fact, in that
case $V(\l)^* \cong V(\l)$ under the mapping which sends
$\l_{n-1}\cdots \l_1\l_0 v_n^*$ to $v_n$.
 We summarize what we have shown for
highest weight modules and the analogous results for lowest
weight and doubly-infinite modules in the following proposition. 
\b 
\proclaim {Proposition 2.35} If $Y$ is a module for $A(\a,\be,\ga)$, 
let $Y^*$ denote the finite dual of $Y$.   
\m
\item {}{(a)} If $V(\l)$ is a highest weight module as in Proposition 2.2, then  
$V(\l)^*$ is an $A(\a,\be,\ga)$-module (which is a highest weight module
isomorphic to $V(\l)$ whenever $\l_n \neq 0$ for all $n$) under the action

$$\aligned & u \cdot v_n^* = \l_n v_{n+1}^*  \\
& d \cdot v_n^* = v_{n-1}^*. \\ 
\endaligned $$
\m 
\item {}{(b)} If $W(\k)$ is a lowest weight module as in Proposition 2.30, then  
$W(\k)^*$ is an $A(\a,\be,\ga)$-module  (which is a lowest weight module
isomorphic to $W(\k)$ whenever $\k_n \neq 0$ for all $n$)  under the action

$$\aligned & u \cdot w_n^* = w_{n-1}^*  \\
& d \cdot v_n^* = \k_n w_{n+1}^*. \\ 
\endaligned $$
\m
\item {}{(c)} If $V(\k,\l)$ is a doubly-infinite module
as in Proposition 2.35, then  
$V(\k,\l)^*$ is a doubly-infinite module under the action

$$\aligned & u \cdot v_n^* = \l_n v_{n+1}^*  \\
& d \cdot v_n^* = v_{n-1}^*. \\ 
\endaligned $$  \endproclaim
 
\b
\head { \S 3. Structural results on down-up algebras} \endhead

\b
\proclaim {Theorem 3.1} (Poincar\'e-Birkhoff-Witt Theorem) Assume $A =
A(\a,\be,\ga)$ is a down-up algebra over $\C$. Then
$\{u^i(du)^jd^k \mid i,j,k = 0,1,\dots \}$ is a basis of $A$.
\endproclaim
\b
{\bf Proof.} It is clear from the
defining relators $f_1 = d(du) - \a (du)d - \be ud^2 - \ga d$
and $f_2 = (du)u -\a u(du) - \be u^2 d - \ga u$ that
the monomials $u^i(du)^jd^k,$ $i,j,k = 0,1, \dots$  are all the irreducible words 
in $A$ and so
span $A$.  
We show that they are linearly independent by
applying the Diamond Lemma (see [Be]).

We can order the words in the
free algebra $\C\la d,u \ra$ by their total degree in $d$ and $u$, and then
lexicographically with $d > u$. Two words in $\C\la d,u \ra$ {\it admit a 
composition} if an end of one is a beginning of the other.  We need to argue
that whenever the largest words $w_1 = d^2 u$ and $w_2 = du^2$ in the relators
$f_1$ and $f_2$ admit a composition say, $w_1 = ab$ and $w_2 = bc$ or 
$w_2 = a'b'$ and $w_1 = b'c'$, then

$$\aligned 
f_1c - af_2 & = a_1 t_1 f_1 x_1 + a_2 t_2 f_2 x_2 \quad \quad \text {or}  \\
f_2c' - a'f_1 & = b_1 y_1 f_1 z_1 + b_2 y_2 f_2 z_2, \\
\endaligned $$

\n where $a_i,b_i \in \C$, $t_i,x_i,y_i,z_i$ are words in $d,u$ (possibly empty),  
$t_if_ix_i < abc$, and $y_if_iz_i < a'b'c'$ for $i = 1,2$.  There is exactly one
admissible composition, $d^2 u^2 = d(du)u$,  of the largest words in the relators, and

$$\aligned 
(d^2 u - \a dud - \be ud^2 - \ga d)& u - d(du^2 - \a udu -\be u^2d - \ga u) \\
& = \be u d^2 u - \be d u^2 d  = 0.\\ 
\endaligned
$$

\n  Therefore by the
Diamond Lemma, the irreducible words form a basis of $A$.    
\qed 
\b
\proclaim {Corollary 3.2} The Gelfand-Kirillov dimension of any down-up
algebra $A(\a,\be,\ga)$ is 3.  \endproclaim
\b
{\bf Proof.}  Let  $A^{(n)} =$ span$_\C\{u^i (du)^j d^k \mid
i+2j+k \leq n\}$.  The spaces $A^{(n)}$ afford a filtration
$(0) \subset A^{(0)} \subset A^{(1)} \subset \dots \subset \cup_n A^{(n)} =
A(\a,\be,\ga)$ of the down-up algebra, and $A^{(m)} A^{(n)} \subseteq A^{(m+n)}$ since the
defining relations replace the words $d^2 u$ and $d u^2$ by words of the
same or lower total degree.   The number of monomials $u^i (du)^j d^k$
with $i+2j+k = \ell$ is $(m+1)(m+1)$ if $\ell = 2m$ and is $(m+1)(m+2)$
if $\ell = 2m+1$.    Thus, dim$A^{(n)}$ is a polynomial in $n$ with
positive coefficients of degree 3, and 
the Gelfand-Kirillov dimension is given by

$$\aligned
\text{GKdim}\big(A(\a,\be,\ga)\big) & = \limsup_{n \rightarrow \infty} \log_n(\dim
A^{(n)})
\\ & = \lim_{n \rightarrow \infty}{\frac {\ln\Big(\dim
A^{(n)}\Big)}{\ln n}} \\ & = 3.
\qed
\\
\endaligned$$
\b
   Recall that the {\bf Jacobson radical} of a ring $R$ is the intersection
of all the annihilators of the simple $R$-modules.  It has many
other characterizations, but we will not need them here.  A ring whose
Jacobson radical is zero is said to be {\bf semiprimitive}.  
\b
\proclaim {Theorem 3.3} If $A(\a,\be,\ga)$
has infinitely many simple Verma modules
$V(\l)$, then the intersection  of the
annihilators of the simple Verma modules is zero.  
Therefore, such a down-up algebra $A(\a,\be,\ga)$ is semiprimitive.
\endproclaim 
\m
{\bf Proof.} Suppose the element $a = \sum_{i,j,k} a_{i,j,k} u^i (du)^j d^k$
annihilates all the simple modules
$V(\l)$.  Assume there is a least
value $\ell$ so that $a_{i,j,\ell} \neq 0$ for some $i,j$.  Let $V(\l)$ be
a simple highest weight module for $A(\a,\be,\ga)$ of weight $\l$.
Applying $a$ to the vector $v_\ell$ in $V(\l)$ we get

$$\aligned
0 = a \cdot v_\ell & = (\l_{\ell-1} \cdots \l_1 \l_0)\sum_{i,j} a_{i,j,\ell} u^i (du)^j
\cdot v_0  \\
& = (\l_{\ell-1} \cdots \l_1 \l_0)\sum_{i,j} a_{i,j,\ell}\l_0^j u^i\cdot v_0 \\
& = (\l_{\ell-1} \cdots \l_1 \l_0)\sum_{i,j} a_{i,j,\ell}\l_0^j v_i. \\
\endaligned $$ 

\n Recall that a necessary and sufficient condition for $V(\l)$ to be simple is
that none of the $\l_n$ is zero.  The elements $v_i$ are linearly independent in
$V(\l)$, so combining those results, we see for each $i$ that

$$0 =\sum_j a_{i,j,\ell} \l_0^j  = \sum_j a_{i,j,\ell} \l^j.$$

\n Because there are infinitely many values of $\l$ that give simple
highest weight modules and only finitely many $\l$ that can be the roots of a nonzero
polynomial, we see that the coefficients $a_{i,j,\ell}$ must all be zero.
Since the Jacobson radical is the intersection
of the annihilators of all the simple $A$-modules, it
too must be zero.   \qed
\b
Let us examine in detail when a down-up algebra $A(\a,\be,\ga)$ has
infinitely many simple Verma modules.  Since the Verma
module $V(\l)$ is simple if and only if $\l_n \neq 0$ for
all $n \geq 1$, the analysis requires the expressions for 
$\l_n$ in Proposition 2.12.
\m
\noindent Case 1. $\a^2 + 4\be \neq 0$.  Here $\l_n = c_1 r_1^n + c_2 r_2^n + x_n$,
(see Proposition 2.12 (i)), 
which can be expressed as a polynomial in $\l$:

$$\l_n = \frac{1}{r_2-r_1}\Big((r_2-\a)r_1^n + (-r_1+\a)r_2^n\Big)\l + b_n,$$

\n where $r_2-\a = 1/2\a -t -\a = -r_1$ and $-r_1+\a = -1/2\a -t + \a = r_2$.
Thus the coefficient of $\l$ is $(r_2-r_1)^{-1}(r_2^{n+1}-r_1^{n+1})$.  If this
is 0, then $r_1^{n+1} = r_2^{n+1}$, which implies $(\a + \tau)^{n+1} = (\a-\tau)^{n+1}$
where $\tau = 2t$ and $\tau^2 = \a^2 + 4 \be$.   Since $\tau^2 = \a^2 + 4\be \neq 0$,  
 $\a+\tau$ and $\a-\tau$ are unequal, and we may suppose one of them --- say $\a+\tau \neq 0$.
Then $1 = \Big((\a+\tau)^{-1} (\a-\tau)\Big)^{n+1}$, so in fact both are nonzero,
and the same holds true if $\a-\tau \neq 0$.    We may
suppose $(\a+\tau)^{-1} (\a-\tau) = e^{i\theta}$ where $\theta = 2\pi k/(n+1)$ for
some $1 \leq k \leq n$.  Then
$\a-\tau = (\a+\tau)e^{i\theta}$ gives $\a(e^{i\theta} -1) = -\tau(e^{i\theta}+1)$. Squaring that
relation shows that $\a^2 = -\be(e^{i\theta}+e^{-i\theta})^2 = -4\be \cos^2(\theta/2)$.  
Consequently, $\l_n$ is a degree 1 polynomial in $\l$ except when
$\a^2 = -4\be \cos^2(\theta/2)$ where $e^{i \theta}$ is an $(n+1)$st root of unity.  
\m
\n Case 2. $\a^2 + 4\be = 0$ and $\a \neq 0$.   In this case $\l_n = c_1s^n + c_2 ns^n + x_n
= \big((n+1)s^n\big)\l + b_n$, where $s = \a/2$.  Thus, $\l_n$ is a degree 1 polynomial in $\l$ for
all
$n$.
\m
\n Case 3. $\a^2+4\be = 0$ and $\a = 0$ (so $\be = 0$).  Here $\l_n = \ga$ for all $n \geq 1$.
\m
\proclaim{Proposition 3.4} If 
\m
\item {}{(a)} $\a^2+4\be \neq 0$ and $\a^2 \neq -4\be \cos^2(\theta/2)$
for some root of unity $e^{i\theta},$  
\s
\item {}{(b)} $\a^2+4\be = 0$ and $\a \neq 0$, or
\s
\item {}{(c)}
 $\a = 0 = \be$ and $\ga \neq 0$,  
\m
\n then $A(\a,\be,\ga)$ has infinitely many simple Verma modules $V(\l)$. 
\endproclaim
\b
{\bf Proof}.  In the cases described in (a) and (b), each $\l_n$ can be regarded as a line
$\l_n = m_n \l + b_n$ with nonzero slope $m_n$.  Each such line intersects the line
$y = 0$ at one value of $\l$.  Avoiding all those countably many values by choosing $\l$
appropriately, we get infinitely many simple Verma modules $V(\l)$.   \qed
\b
\proclaim {Proposition 3.5} A down-up algebra $A = A(\a,\be,\ga)$ is
$\Z$-graded, $A = \oplus_n A_n$, where $A_n =$ span$_{\C}\{u^{i} (du)^j d^k \mid i,j,k =
0,1,\dots$, and $i-k = n\}$, and 
$A_0$
is a commutative subalgebra. \endproclaim
\b
{\bf Proof.} The free associative algebra over $\C$ generated by $d$ and $u$
can graded by assigning  $\deg(d) = -1$ and $\deg(u) =1$ and extending this
by using $\deg(ab) = \deg(a) + \deg(b)$.  The relations
$d^2 u = \a dud + \be ud^2 + \ga d$ and $du^2 = \a udu + \be u^2 d + \ga u$
are homogeneous, so the down-up algebra inherits the grading.  

To prove that $A_0$
is commutative, we first show that $u^i (du)^j d^i$ commutes with $du$ and $ud$ for all
$i,j$.  This is true if $i+j \leq 1$ by (2.11).  Suppose we know the result for all $i+j < N$,
and let $i+j = N \geq 2$.   Consider $(u^i (du)^j d^i) (ud)$.  Now if
$i = 0$ this equals $(ud)(u^i (du)^j d^i)$ since $ud$ and $du$ commute. If $i = 1$,
then $u(du)^jd = (ud)^{j+1}$, which commutes with $ud$.  So we may assume that $i \geq 2$.
Then

$$\aligned
(u^i (du)^j d^i)(ud) & = u^i(du)^j d^{i-2} (\a dud + \be ud^2 + \ga d)d  \\
& = \a u^2 (u^{i-2}(du)^jd^{i-2})(du)d^2 + \\
& \hskip .5 truein  \be u^2(u^{i-2}(du)^jd^{i-2})(ud)d^2
+ \ga u^i (du)^j d^i \quad \quad \text{while} \\
\endaligned$$

$$\aligned
(ud)(u^i (du)^j d^i) & = u(\a udu + \be u^2d + \ga u)(u^{i-2}(du)^j d^i)  \\
& = \a u^2 (du) (u^{i-2}(du)^jd^{i-2})d^2 + \\
& \hskip .5 truein \be u^2(ud)(u^{i-2}(du)^jd^{i-2})d^2
+ \ga u^i (du)^j d^i. \\
\endaligned$$ 

\n These two expressions are equal by our induction hypothesis.  Analogously,
we consider $(u^i (du)^j d^i) (du)$, where we can suppose that $i \geq 1$. Then

$$\aligned
(u^i (du)^j d^i)(du) & = u^i(du)^j d^{i-1} (\a dud + \be ud^2 + \ga d) \\
& = \a (u^i (du)^jd^i)(ud) + \\
& \hskip .5 truein  \be u(u^{i-1}(du)^jd^{i-1})(ud)d
+ \ga u^i (du)^j d^i \quad \quad \text{which equals} \\
\endaligned$$

$$\aligned
(du)(u^i (du)^j d^i) & = (\a udu + \be u^2d + \ga u)(u^{i-1}(du)^j d^i)  \\
& = \a (ud)(u^i(du)^jd^i) + \\
& \hskip .5 truein \be u(ud)(u^{i-1}(du)^jd^{i-1})d
+ \ga u^i (du)^j d^i
 \\
\endaligned$$ 

\n by what we have already established. 

The proof of the proposition will be complete once we show that $u^i(du)^jd^i$
commutes with $u^k(du)^\ell d^k$ for all $i,j,k,\ell$.  This is true if
$i+j \leq 1$ or if $i$ or $k$ is $\leq 1$ by our previous considerations, so we assume it is
true when
$i+j+k+\ell < N$ and prove it for $i+j+k+\ell = N$.  We suppose $i,k \geq 2$. 
Now

$$\aligned (u^i(du)^jd^i)(u^k(du)^\ell d^k) & = u^i(du)^jd^{i-2}
(\a dud + \be ud^2 + \ga d)u^{k-1}(du)^\ell d^k \\
& = \a u(u^{i-1}(du)^jd^{i-1})(ud)(u^{k-1}(du)^\ell d^{k-1})d + \\
& \hskip .5 truein  \be u^2(u^{i-2}(du)^jd^{i-2})(ud)(du)(u^{k-2}(du)^\ell d^{k-2})d^2 + \\
& \hskip .5 truein  u(u^{i-1}(du)^jd^{i-1})(u^{k-1}(du)^\ell d^{k-1})d, \\
\endaligned $$

\n which can be seen to equal the same expression with the roles of $i,j$ and $k,\ell$
reversed by induction.   \qed

\b
\subhead The center of a down-up algebra \endsubhead
\b
Let $Z(A) = \{z \in A \mid za = a z$ for all $a \in A\}$ denote the
center of a down-up algebra $A = A(\a,\be,\ga)$.   
\b
\proclaim {Proposition 3.6} Suppose $A = \sum_{j \in \Z} A_j$
is the $\Z$-grading of $A =A(\a,\be,\ga)$ as in Proposition 3.5, and
assume that $A$ has infinitely many simple Verma modules.  Then: 
\item {}{(a)} The center $Z(A)$ of $A$ is contained in $A_0$ = span$_\C\{u^i(du)^jd^i \mid i,j =
0,1, \dots \}$. 
\m
\item{}{(b)}  
If $z \in Z(A)$ and $V(\l)$ is any Verma module of
$A$, then $z$ acts as a scalar, say $\chi_\l(z)$, on  
 $V(\l)$.
The mapping $\chi_\l : Z(A) \rightarrow \C$ is an algebra homomorphism.  
\m
\item {} {(c)} A scalar  $\pi \in \C$ is {\bf linked} to $\l$ if $\pi = \l_{n+1}$
for some $n \geq 0$ with $\l_n = 0$, where the sequence
$\l_1, \l_2, \dots$  is 
constructed using the recurrence relation in (2.1) starting with  
$\l_{-1} = 0$ and $\l_0 = \l$.  If $\pi$ is linked to $\l$, then
$\chi_\pi = \chi_\l$.  
\endproclaim
\b
{\bf Proof.} Consider
the sequence $\l = \l_0,\l_1,\l_2, \dots$ from (2.1),
and let  $\Lambda(\l) = \{n \in \Z_{\geq 0} \mid \l_n = 0\}$.  
Suppose $v_0, v_1, \dots$ is the canonical basis for $V(\l)$ as
in Proposition 2.2, and let $z \in Z(A)$.   From 
$dz \cdot v_0 = zd \cdot v_0 = 0$ 
we see that $z \cdot v_0 = \chi_\l(z) v_0 + \sum_{n \in \Lambda(\l), n > 0} 
a_n v_n$ for some scalars $a_n$ and $\chi_\l(z)$.  Therefore,
when $\l \neq 0$ and $V(\l)$ is simple (i.e. when $\Lambda(\l) = \emptyset$)
then $z \cdot v_0 = \chi_\l(z) v_0$.  But  $z a \cdot v_0 = a z \cdot v_0 = 
\chi_\l(z) a \cdot v_0$ for all $a \in A$, and since $V(\l) = A v_0$,
this shows that $z$ acts as the scalar $\chi_\l(z)$ on such Verma modules.
Now write $z$
as a sum of its homogeneous components, 
$z =\sum_j z_j$, $z \in A_j$.   Then  
 since $A_j v_n \subseteq \C v_{j+n}$
for all $j,n$, and since $z v_n = \chi_\l(z) v_n$ for all $n$ when $V(\l)$
is simple, it must be
that $z_j$ for $j \neq 0$ lies in the annihilator of $V(\l)$.  But
the intersection of the annihilators of the simple $V(\l)$ modules
is zero by Theorem 3.3.  Thus, $z = z_0 \in A_0$, as claimed in (a).
Once we know that, then it follows for any Verma module $V(\l)$ that
$z \cdot v_0 = \chi_\l(z)v_0$ for some scalar $\chi_\l(z) \in \C$.  As before
using $V(\l) = A v_0$, we see $z$ acts as $\chi_\l(z)$ on all of $V(\l)$.
As $\chi_\l(zz')v_0 = z z' \cdot v_0 = z \cdot (z' v_0) = \chi_\l(z)\chi_\l(z')v_0$,
it is clear $\chi_\l$ is an algebra homomorphism.  Part (c) is apparent
also, since when $\l_n = 0$, then span$_\C\{v_j \mid j \geq n+1\}$ is
a submodule of $V(\l)$ isomorphic to $V(\l_{n+1})$.    \qed

\b
\head { \S 4. Category $\Cal O$ modules} \endhead
\b
The representation theory of finite-dimensional complex semisimple 
Lie algebras centers on a certain important collection of weight modules,
the so-called category $\Cal O$ modules, whose study was initiated 
by Bernstein, Gelfand, and Gelfand [BGG].  In this section
we define an analogous
category of modules for each down-up algebra $A = A(\a,\be,\ga)$.
We begin with a general result about weight vectors and then
specialize to the modules in $\Cal O$ and derive some
results about them.  
\b
\proclaim {Proposition 4.1} Suppose $M$ is an $A$-module and $m \in M$ is a vector  
of weight $\nu = (\nu',\nu'')$.  
\s
\item{}{(a)} Then $u\cdot m$ is a vector of weight $\mu = (\mu',\mu'')$
for $\mu' = \a \nu' + \beta \nu'' + \ga$ and $\mu'' = \nu'$.  
\m
\item{}{(b)} If $\beta\neq 0$, then $d\cdot m$ is a vector of weight $\delta =
(\delta',\delta'')$ where $\delta' = \nu''$ and $\delta'' = \beta^{-1}(\nu' - \a \nu'' -
\ga)$.  \endproclaim
 
\b
{\bf Proof.} The assertions are evident from the following calculations:

$$
\gather\tag 4.2 \\
(du)u \cdot m  = \a u(du) \cdot m + \be u(ud) \cdot m + \ga u \cdot m = (\a \nu' + \be
\nu'' +
\ga) u \cdot m
\\ (ud)u \cdot m   = u(du) \cdot m = \nu' u\cdot m \\
(du)d \cdot m   = d(ud) \cdot m = \nu'' d\cdot m \\
(ud)d\cdot m  = \beta^{-1}\big(d(du)\cdot m - \a d(ud) \cdot m - \ga d \cdot m \big) = 
\beta^{-1}\big(\nu' - \a \nu'' - \ga \big) d\cdot m.  \qed  \\  
\endgather $$ 
\b 
\n (4.3)~~The category $\Cal O$ consists of all $A$-modules $M$ satisfying the
following conditions:
\m
\item {}{(a)} $M$ is a weight module
relative to $\h =$ span$_\C\{du,ud\}$, i.e.  $M = \bigoplus_{\nu} M_\nu$
where $M_\nu = \{m \in M \mid h \cdot m = \nu(h) m$ for all $h \in \h \}$;
\s
\item {}{(b)} $d$ acts locally nilpotently on $M$, so that for each
$m \in M$, $d^n \cdot m = 0$ for some $n$.  
\s
\item {}{(c)} $M$ is a finitely generated $A$-module. 
\m
The category $\Cal O$ is closed under taking submodules
and quotients.  It contains all the Verma modules $V(\l)$
and their simple quotient modules $L(\l)
= V(\l)/M(\l)$ (and in the case
that  $\gamma = 0$ and $\l = 0$, the one-dimensional quotients
$L(0,\xi) \eqdef  V(0)/N^{(x-\xi)}$ (see Corollary 2.28)).   
Modules in the category $\Cal O$ enjoy some of the same basic properties
as the category $\Cal O$ modules for complex semisimple Lie algebras
and for the algebras studied by Smith [Sm].
\b
\n (4.4) Suppose $\w = (\w',\w'')$, where
$\w',\w'' \in \C$.  {\it Let $\Omega(\w)$ denote the set
weights defined inductively starting with 
$\w_0 = (\w_0',\w_0'')$, where $\w_0' = \w'$ and $\w_0'' = \w''$,
and proceeding for $n \geq 1$ by setting $\w_n = (\w_n',\w_n'')$ where}

$$\aligned \w_n' & = \a \w_{n-1}' + \be \w_{n-1}'' + \ga \\
\w_n'' & = \w_{n-1}'. \\
\endaligned$$

\n Note that $\w_n' = \a \w_{n-1}' + \be \w_{n-2}' + \ga$
for $n \geq 2$.  When $\omega = (\l,0)$
then $\w_n = (\l_n,\l_{n-1})$,
where $\l_n$ is as in Proposition 2.1.
Thus $\Omega(\w)$ is just the set of weights 
of the Verma module $V(\l)$ in the particular instance that $\w = (\l,0)$.  
\b
\proclaim {Lemma 4.5} Suppose $\beta \neq 0$.  Assume $M$ is an object in the category
$\Cal O$, and let $\Omega(M) = \{ \nu \in \h^* \mid M_\nu \neq (0)\}$ be the set
of weights of $M$.   Then there exists a
finite set of weights $\w^1, \dots, \w^r$ so that $\Omega(M) \subseteq 
\cup_i \Omega(\w^i)$.   \endproclaim 
\b
{\bf Proof.} Because $M$ is finitely generated, there exists a finite set
of weight vectors $m_1, \dots, m_\ell$ generating $M$.  Then $\sum_j \C[d]m_j$ 
is finite-dimensional, and moreover it is $\h$-invariant for
$\h =$ span$_\C\{du,ud\}$ by Proposition 4.1.  So we
can choose
a basis $y_1, \dots, y_r$ for it such that 
$y_i$ has some weight $\w^i$ for each $i$.    Then by Theorem 3.1, 

$$M = \sum_j A m_j = \sum_j \C[u] \C[du] \C [d] m_j = \sum_i \C[u] \C[du] y_i
= \sum_i \C[u] y_i. $$

\n The vectors in $\C[u]y_i$ have weights in $\Omega(\w^i)$ by Proposition 4.1, 
so the relation $\Omega(M) \subseteq 
\cup_i \Omega(\w^i)$ is apparent.   \qed 
\b
\proclaim {Lemma 4.6} Suppose $\beta \neq 0$. If $M$ is a simple object in the category
$\Cal O$, then $M \cong L(\l)$ for some $\l$, or 
$\gamma = 0$ and $M \cong V(0)/N^{(x-\xi)} = L(0,\xi)$ for some $\xi \in \C$
(as in Corollary 2.28).     
\endproclaim
\b
{\bf Proof.}  It follows from the local nilpotence of $d$ 
and from Proposition 4.1 that $M$ contains some highest weight vector
$v_0$.  The simplicity of $M$ forces $M = Av_0$ so that $M$ is a highest
weight module.  As highest weight modules are quotients of Verma modules,
$M$ is as asserted.  \qed 
\m 
\b
\head { \S 5. Category $\Cal O'$ modules} \endhead
\b
In this section we introduce and investigate a more general category of modules
for the down-up algebra $A = A(\a,\be,\gamma)$.  {\bf We require throughout that
$\be \neq 0$.}  
\b
\n (5.1)~~The category $\Cal O'$ consists of all $A$-modules $M$ satisfying the
following conditions:
\m
\item {}{(a)} $M$ is a weight module
relative to $\h =$ span$_\C\{du,ud\}$ 
\s
\item {}{(b)} $\C[d]m$ is finite-dimensional for each $m \in M$.    
\s
\item {}{(c)} $M$ is a finitely generated $A$-module. 
\m   
The modules in $\Cal O$ clearly belong to $\Cal O'$, but we will show later that $\Cal O'$
is larger than
$\Cal O$ by determining the simple modules in $\Cal O'$.  
Two mappings on weights play a prominent role in these investigations. 
\b
Recall (Proposition 4.1) that if a vector $m$ in an $A$-module has weight $\nu =
(\nu',\nu'')$ relative to
$\h$, then $d\cdot m$ has weight 

$$\aligned  \delta(\nu) & = (\delta(\nu)', \delta(\nu)'') \quad \quad \text {where} \\
\delta(\nu)' & = \nu'' \quad \quad \text {and}\quad \quad 
\delta(\nu)'' = \be^{-1}(\nu' - \a \nu'' - \ga), \\
\endaligned \tag 5.2$$

\n and $u\cdot m$ has weight 

$$\aligned  \mu(\nu) & = (\mu(\nu)', \mu(\nu)'') \quad \quad \text {where} \\
\mu(\nu)' & = \a \nu'+ \be \nu'' + \ga \quad \quad \text {and} \quad \quad
\mu(\nu)'' = \nu'. \\
\endaligned \tag 5.3$$

\n An easy direct computation shows that $\delta(\mu(\nu)) = \nu$ 
and $\mu(\delta(\nu)) = \nu$,
(or this can be readily seen from the fact that $(du) \cdot m = \nu'm$ and $(ud)
\cdot m = \nu''
m$).  Observe that the
set of weights $\Omega(\w)$
in (4.4) alternately can be
described as  

$$\Omega(\w) = \{\mu^k(\w) \mid k = 0,1, \dots\}. \tag 5.4$$  
\b
\n The set of weights of the lowest weight module $W(\k)$
is just $\{\delta^k(\w) \mid  k = 0,1, \dots\}$, where
$\w = (0,\k)$.  
\b
\proclaim{Theorem 5.5} Suppose $F$ is a set of weights such that
$\delta(\w), \mu(\w) \in F$ whenever $\w \in F$. 
Suppose $\rho \in \C$ is nonzero, and let $N(F,\rho)$ be the $\C$-vector
space with basis $\{v_\w \mid \w \in F\}$.   \m
\item{}{(a)} Define

$$\aligned d \cdot v_\w & = \rho v_{\delta(\w)} \\
u \cdot v_\w & = \rho^{-1} \mu(\w)'' v_{\mu(\w)}. \\
\endaligned \tag 5.6 $$ 

\n \item {}{}Then this action extends to give an $A(\a,\be,\ga)$-module action on
$N(F,\rho)$. 
\m
\item{}{(b)} If $F$ is generated by any weight $\nu = (\nu',\nu'')  \in F$ under the action
of $\delta$ and $\mu$, and if $\nu' \neq 0$ for any $\nu  \in F$, then $N(F,\rho)$ is a simple
$A(\a,\be,\ga)$-module.  
\endproclaim 
\b
{\bf Proof.}  Let $N(F,\rho)$ be as in the statement of the theorem,
and define an action of the free associative algebra $\C\la d,u\ra$ on
$N(F,\rho)$ using (5.6).  Then

$$
\aligned (du) \cdot v_\w & = \rho^{-1}\mu(\w)'' d \cdot v_{\mu(\w)} = 
\mu(\w)'' v_\w = \w'
v_\w 
\\
(ud) \cdot v_\w & = \rho u \cdot v_{\delta(\w)} = \mu(\delta(\w))''v_\w = \w'' v_\w. \\
\endaligned \tag 5.7$$

\n Utilizing (5.7) together with (5.2) and (5.3), it is easy to see that

$$ \aligned 
& \Big(d^2 u - \a dud - \be u d^2 - \ga d\Big) \cdot v_\w =
\rho\Big(\w'-\a \w'' - \be \delta(\w)'' - \ga \Big)v_{\delta(\w)} = 0 \\
& \Big(du^2  - \a udu - \be u^2 d - \ga u\Big) \cdot v_\w =
\rho^{-1}
\mu(\w)''\Big(\mu(\w)'-\a \w' - \be \w'' - \ga \Big)v_{\mu(\w)} = 0. \\
\endaligned $$

\n Thus, there is an induced action making
$N(F,\rho)$ into an $A(\a,\be,\ga)$-module.  
\b 
Now suppose $F$ is generated by any weight $\nu \in F$ under the maps $\delta$ and
$\mu$ defined in (5.2) and (5.3). It follows
from (5.7) that $N(F,\rho)$ is a weight module. 
Let $M$ be a nonzero submodule of $N(F,\rho)$, and 
assume  $v = \sum_\w a_\w v_\w$ is a nonzero vector in $M$.  Since any submodule of a
weight module is a weight module, 
 $v_\w \in M$ if $a_\w \neq 0$.  
 Now whenever $v_\w \in M$ for some
$\w$, then
$d \cdot v_\w = \rho v_{\delta(\w)}$ and $u \cdot v_\w = \rho^{-1}\mu(\w)'' v_{\mu(\w)}$
both belong to $M$.  
Since $\mu(\w)'' = \w' \neq 0$ for any $\w \in F$
by assumption, $v_{\mu(\w)} \in F$.   
Because $F$ is generated by $\w$ under $\delta$ and
$\mu$,  $M$ contains  $v_\tau$ for all $\tau
\in F$.   As a consequence, $M = N(F,\rho)$ and 
$N(F,\rho)$ is a simple
$A(\a,\be,\ga)$-module.  \qed
\b
Suppose now that $M$ is a simple module in the category $\Cal O'$. 
If $M$ contains a highest weight vector, then $M$ is as in Lemma 4.6.
So we may assume henceforth that $M$ has no highest weight vectors.  
Let $m$ be a nonzero weight vector in $M$ and consider the finite-dimensional
space $\C[d]m$.  If $d^k\cdot m = 0$ for some $k$, then 
$d^{k-1}\cdot m$ is a highest weight vector in $M$ contrary
to assumption.  It must be that $d$ acts
nonnilpotently on $\C[d]m$.  Let  $v$ be an
eigenvector of $d$ in $\C[d]m$ corresponding to
nonzero eigenvalue, say $d\cdot v = \rho v$, $\rho \neq 0$.  Write
$v = \sum_\w v_\w$ where $v_\w$ belongs to the weight space $M_\w$.  

Since  $d \cdot v = \rho v$,  it must be that $d \cdot v$ has the same number of
weight components as $v$.  Equating corresponding components on both
sides of the equation gives

$$d \cdot v_\w = \rho v_{\delta(\w)}. \tag 5.8$$

\n Thus if $S$ is the $\C$-span of the weight components
$v_\w$ of $v$, then $d S = S$.  Now 
$\rho u \cdot v = (ud) \cdot v = \sum_\w \w'' v_\w$,  which implies that

$$u \cdot v_\w = \rho^{-1} \mu(\w)'' v_{\mu(\w)}. \tag 5.9$$

\n  Thus, $S$ is an $A(\a,\be,\ga)$-submodule of $M$.  By
simplicity,  $S = M$.  

Since there are
only finitely many nonzero weight component summands, $\delta^\ell(\w) =
\delta^k(\w)$ for
some $\ell > k$.  Applying $\mu^k$ to both sides, we see $\delta^{(\ell-k)}(\w) = \w$. 
Thus,  there is some smallest value $p$ so that
$\delta^p(\w) = \w$ and  $v_\w, v_{\delta(\w)}, \dots, v_{\delta^{p-1}(\w)}$
are all nonzero.  Their $\C$-span $S'$ 
is invariant under $d$.  Using (5.9) and the fact that $\mu(\delta^i(\w)) =
\delta^{i-1}(\w)$ for $i \geq 1$ and $\mu(\w) = \mu(\delta^{p}(\w)) = \delta^{p-1}(\w)$,
we see that $u S' \subseteq S'$ also.  Thus, $M = S'$ by simplicity.  
The set $F = \{\w,\delta(\w), \dots, \delta^{p-1}(\w)\}$ is $\delta$
and $\mu$-invariant.
Any weight in it generates the whole set under $\delta$ and $\mu$.
Now suppose $\nu$ is a weight of $F$
and $\nu' = \mu(\nu)'' = 0$.  
Setting $v_n = \rho^n v_{\delta^n(\nu)}$ for $n = 0,1,\dots, p-1$,
we have $u \cdot v_0 = 0$, ~ $d \cdot v_n = v_{n+1}$ for $0 \leq n \leq p-1$, 
and $d \cdot v_n = \rho^p v_0$.   
In this case $M$ is a finite-dimensional lowest weight module
-- a quotient of $W(\k)$ for $\k = \nu''$. 
The sole remaining possibility
is that $\nu' \neq 0$ for any $\nu \in F$. In this situation 
$M = N(F,\rho)$, (and  $N(F,\rho)$ is in fact a simple module - compare
Theorem 5.5 (b)).    
\b
Let us summarize what we have shown:
\b
\proclaim {Theorem 5.10}  Assume $M$ is a simple module in the category $\Cal O'$.
Then there are three possibilities:
\m
\item{}{(a)} $M$ is a highest weight module, that is, $M$ is isomorphic to
 $L(\l)$ for some $\l$ or to $L(0,\xi)$ for some $\xi \in \C$ (when $\ga = 0$).
\m
\item{}{(b)} $M$ is a finite-dimensional lowest weight module with
weights $\nu, \delta(\nu), \dots, \delta^{p-1}(\nu)$ such that
 $\delta^p(\nu) = \nu$.
\m
\item{}{(c)} $M$ is isomorphic to $N(F,\rho)$ for some $\rho \neq 0$ 
and some finite set $F
= \{\nu, \delta(\nu), \dots, \delta^{p-1}(\nu)\}$ such that
 $\delta^p(\nu) = \nu$. 
\endproclaim
\b
We conclude this section by analyzing the special case that the
set $F$ is cyclically generated by $\w = (\l,\k)$ under $\delta$ and $\mu$,
and the weights in $$F = \{ \dots, \mu^2(\w),\mu^1(\w),\w,\delta(\w),
\delta^2(\w), \dots \}$$ are all distinct.  Observe that if
we define $\l_n$ by first setting $\l_0 = \l$ and $\l_{-1} = \k$ and proceeding
by

$$\aligned
\l_n & = \a \l_{n-1} + \be \l_{n-2} + \ga  \quad \quad \text {for} \; n \geq 1 \\
\l_{-n}&  = \be^{-1}(\l_{-n+2} -\a \l_{-n+1} -\ga) \quad \quad \text {for} \; n \geq 2, \\
\endaligned $$

\n then

$$\aligned (\l_n,\l_{n-1}) & = \mu^n(\w)  \quad \quad \text {for} \; n \geq 0 \\
(\l_{-n},\l_{-n-1}) & = \delta^n(\w)  \quad \quad \text {for} \; n \geq 1. \\
\endaligned $$

\n These are just the weights of the doubly-infinite module $V(\k,\l)$ in
Proposition 2.33.  It is easy to see
using Proposition 2.35 that there is an isomorphism from the dual module $V(\k,\l)^*$
to $N(F,1)$ taking  $v_n^* \mapsto v_{\mu^n(\w)}$ for $n \geq 0$
and $v_{-n}^* \mapsto v_{\delta^n(\w)}$ for $n \geq 1$. Theorem 5.5 shows
that if $\nu' = \l_n  \neq 0$ for any $\nu = (\l_n,\l_{n-1})$, $n \in \Z$,
then $N(F,1) \cong V(\k,\l)^*$ is simple.   This forces
$V(\k,\l)$ to be simple as well, for the annihilator in $V(\k,\l)^*$ 
of any proper submodule of $V(\k,\l)$ would be a proper submodule of $V(\k,\l)^*$.  
Consequently, an analogue of the result in Proposition 2.4 holds for the
doubly-infinite module $V(\k,\l)$:
\b
\proclaim {Proposition 5.11} Suppose $\l_0 = \l$ and $\l_{-1} = \k$ and
define the sequence $\l_n$ as in Proposition 2.33.  If 
$\l_n \neq 0$ for any $n$, then the doubly-infinite module
$V(\k,\l)$ is simple.  \endproclaim 
\m 
\b
\head {\S 6. One-dimensional $A(\a,\be,\ga)$-modules and isomorphisms} \endhead
\b
Consider a one-dimensional space  $X = \C x$ 
with

$$d \cdot x = a x \quad \quad \text {and} \quad \quad u \cdot x = b x. $$

\n
Then this extends to a module for $A = A(\a,\be,\ga)$ if and only if 

$$\gather
0 = (d^2 u -\a dud - \be u d^2 - \ga d) \cdot x = a\big((1-\a-\be)ab - \ga\big)x \\
0 = (du^2 -\a udu - \be u^2 d - \ga u) \cdot x = b\big((1-\a-\be)ab - \ga\big)x. 
\\
\endgather$$

If $a = 0 = b$, then $X = V(0)/M(0)$. 
We may suppose that not both $a,b$ are 0, so that $X$ is an $A$-module if and only if
$(1-\a-\be)ab = \ga$ holds.
When $\ga = 0$ and $\a + \be = 1$, then
any choice of $a$ and $b$ gives $X$ the structure of an $A$-module, so
the one-dimensional $A$-modules are parametrized by $\C^2$ in this
instance.    
When $\ga = 0$ and $\a+\be \neq 1$, then either $a = 0$ or $b = 0$ must
hold.  The one-dimensional $A$-modules are parametrized by
$\{(a,0) \mid a \in \C\} \cup \{(0,b) \mid b \in \C\}$ in this case. 
When $\ga \neq 0$, we may suppose that $\ga = 1$.  
>From $(1-\a-\be)ab = \ga = 1$, we see that  $b = a^{-1}(1-\a-\be)^{-1}$.  Thus, 
the one-dimensional $A$-modules are parametrized by $\C$ for $\ga \neq 0$ and $\a+\be \neq
1$,
and by $a = 0 = b$ when $\ga \neq 0$ and $\a+\be = 1$. 

Now $(du) \cdot x = ab x = (ud) \cdot x$, so that $x$ has weight $\nu = (ab,ab)$
for all one-dimensional $A$-modules $X$.  When $\be \neq 0$, then
$\delta(\nu) = (ab, \be^{-1}(ab - \a ab - 1)$, which equals $(ab,ab)$ when $\ga = 1$.
Thus, in the case $\be \neq 0, \ga \neq 0$ the one-dimensional module
$X$ is isomorphic to $N(F,a)$ where $F = \{\nu = (ab,ab)\}$ and $\delta(\nu) = \nu$
(compare Theorem 5.10).   To summarize we have:
\b
\proclaim {Theorem 6.1} The one-dimensional modules for the down-up algebra
$A(\a,\be,\ga)$ are parametrized by
\m
\item {}{(a)} $\C^2$ when $\ga = 0$ and $\a + \be = 1$;
\s
\item {}{(b)} $\{(a,0) \mid a \in \C\} \cup \{(0,b) \mid b \in \C\}$ when $\ga = 0$
and $\a + \be \neq 1$;
\s
\item {}{(c)} $\C$ when $\ga \neq 0$ and $\a + \be \neq 1;$ 
\s
\item {}{(d)} $(0,0)$ when $\ga \neq 0$ and $\a + \be = 1$.  
\endproclaim
\b
\proclaim {Corollary 6.2} If the down-up algebras $A(\a,\be,\ga)$ 
and $A(\a',\be',\ga')$ are isomorphic, then either both $\a+\be$ and $\a'+\be'$
are equal to 1 or both are different from 1, and either both $\ga$ and $\ga'$
are 0 or both are different from 0.  \endproclaim
\b
{\bf Proof.} The result is a simple consequence of Theorem 6.1 and the fact that if
$\phi: A(\a,\be,\ga) \rightarrow 
A(\a',\be',\ga')$ is an isomorphism, then any one-dimensional
module $X = \C x$ for $A(\a',\be',\ga')$ becomes one for $A(\a,\be,\ga)$ 
by defining  $t \cdot x = \phi(t) \cdot x$ for all $t \in A(\a,\be,\ga)$.   \qed

\m
\b
{\bf Open Problems.} \quad In conclusion we mention some noteworthy open
questions concerning down-up algebras which are not addressed here. Problems
(a)-(f) were posed in [B], and since the original version of this
paper appeared, several of these problems have been solved.     
\b
\item {(a)} {\it Determine when the algebra $A(\a,\be,\ga)$ is Noetherian.}  
\quad   
Kirkman, Musson, and Passman [KMP] have recently proven 
that $A(\a,\be,\ga)$ is left and right Noetherian if and only if $\be \neq 0$ if
and only if  
$du$ and $ud$ are algebraically independent.    
\m
\item {(b)} {\it Find conditions on $\a,\be,\ga$ for $A(\a,\be,\ga)$
to be a domain.} \quad   When $\be = 0$, then $d(du - \a ud -\ga 1) = 0$
so that  $A(\a,\be,\ga)$ has zero divisors for any choice of $\a,\ga \in \C$.
It has been shown in [KMP] that
$A(\a,\be,\ga)$ is a domain
if and only if $\be \neq 0$ (i.e if and only if
the conditions in (a) hold.)  It has been proved
independently by Kulkarni [K2] (using methods
from hyperbolic rings) that $\be \neq 0$ is a necessary and sufficient condition
for $A(\a,\be,\ga)$ to be a domain.   
\m
\item {(c)} {\it Compute the center of  $A(\a,\be,\ga)$.} \quad  
This has been solved
in a very recent preprint by Zhao [Z] and announced by Kulkarni [K2].
\m
\item {(d)} {\it  Determine when $A(\a,\be,\ga)$ is a Hopf algebra.} \quad   The
mapping $S(d) = -d$, $S(u) = -u$, which is an algebra antiautomorphism (in fact,
the antipode)
in the enveloping algebra case, is an algebra antiautomorphism of $A(\a,\be,\ga)$
if and only if $\be = -1$ or $\be = 1$ and $\a = 0 = \ga$.
\m 
\item {(e)} {\it Relate Kulkarni's presentation of the maximal left ideals in conformal
$\fsl$ algebras to the simple modules in category $\Cal O$ and category $\Cal O'$.}
\quad  
The approach in [K1] using noncommutative algebraic geometry is quite different from the one
adopted here and so is the description of the simple modules. 
Our primary focus here was on determining
explicit information about the
representations of down-up algebras that could be used in the study of posets
in the spirit of the work of Stanley [St1] and Terwilliger [T].   
\m
\item {(f)} {\it Study the homogenization $A[t]$ of the down-up algebra $A = A(\a,\be,\ga)$,
which is the graded algebra generated by $d,u,t$ subject to the relations} 

$$\gathered d^2 u = \a dud + \be ud^2 + \ga dt^2, \quad  du^2 =
\a udu + \be u^2 d + \ga ut^2, \\
dt = td, \quad ut = tu.\\ \endgathered $$

\item{} Homogenized $\fsl$ is a positively
graded Noetherian domain and a maximal order,
which is Auslander-regular of dimension 4 and satisfies the Cohen-Macaulay
property.  Le Bruyn and
Smith [LS] have determined the point, line, and plane modules of homogenized $\fsl$
and shown the line modules are homogenizations of the Verma modules.
\m 
\item{(g)} {\it Investigate indecomposable and projective modules
for down-up algebras.}
\m
\item{(h)} {\it Determine when two down-up algebras are isomorphic.} \quad  
Some information on this problem is derived in Section 6 as
a consequence of the determination of the 1-dimensional modules, but we do
not address the general problem here.  
\m
\item{(i)} {\it Investigate prime and primitive ideals of down-up algebras.}

\vskip 5 mm
 
\address 
\newline 
Georgia Benkart,
Department of Mathematics, University of Wisconsin, Madison,
Wisconsin 53706-1388   \newline
benkart\@math.wisc.edu
 \newline \newline
Tom Roby,
Department of Mathematics and Computer Science, California State 
University, Hayward, California 94542-3092
\newline
troby\@mcs.csuhayward.edu
\endaddress 

\enddocument